\documentclass[12pt,thmsa]{article}
\usepackage{amsfonts, amsmath}

\newtheorem{theorem}{Theorem}
\newtheorem{corollary}{Corollary}
\newtheorem{proposition}{Proposition}
\newtheorem{lemma}{Lemma}

\newcommand{\p}{\Bbb{P}}

\newcommand{\e}{\Bbb{E}}

\newcommand{\ind}{\mbox{\rm 1\hspace{-0.04in}I}}

\newcommand{\R}{\mbox{\rm I\hspace{-0.02in}R}}

\newcommand{\ed}{\stackrel{(d)}{=}}
\newcommand{\ud}{\mathrm{d}}
\newcommand{\eqdef}{\stackrel{\mbox{\tiny$($def$)$}}{=}}

\def\QED{\hfill\vrule height 1.5ex width 1.4ex depth -.1ex \vskip20pt}

\begin{document}
\hspace*{-0.5in} {\footnotesize This version April 5, 2006.}
\vspace*{0.9in}
\begin{center}
{\LARGE  On the future infimum of positive self-similar
Markov processes.\vspace*{0.4in}}\\
{\large J.C. Pardo \footnote{Research suported by a grant from
CONACYT (Mexico).}\vspace*{0.2in}}
\end{center}
\noindent $^1$  {\footnotesize Laboratoire de Probabilit\'es et
Mod\`eles Al\'eatoires, Universit\'e Pierre et Marie Curie,
4, Place Jussieu - 75252 {\sc Paris Cedex 05.} E-mail: pardomil@ccr.jussieu.fr}\\

\noindent {\it Abstract} {\footnotesize We establish integral
tests and laws of the iterated logarithm for the upper envelope of
the future infimum of positive self-similar Markov processes and
for increasing self-similar Markov processes at 0 and $+\infty$.
Our proofs are based on the Lamperti representation and time
reversal arguments due to  Chaumont and Pardo \cite{ChP}. These
results extend laws of the iterated
logarithm for the future infimum of Bessel processes due to Khoshnevisan et al. \cite{Kh}.}\\


\noindent
{\it A.M.S. Classification}: {\footnotesize 60 G 18, 60 G 17, 60 G 51, 60 F 15.}\\

\section{Introduction and main results}
An $\R_{+}$-valued Markov process $X=(X_{t}, t\geq 0)$ with
c\`adl\`ag paths is a self-similar process if for every $k>0$ and
every initial state $x\geq 0$ it satisfies the scaling property,
i.e., for some $\alpha>0$
\begin{equation*}
\textrm{the law of } (kX_{k^{-\alpha}t},t\geq 0) \textrm{ under
}\p_{x}\textrm{ is } \p_{kx},
\end{equation*}
where $\p_{x}$ denotes the law of the process $X$ starting from $x\geq 0$.\\
We will refer to positive self-similar Markov processes as PSSMP. We
will also denote by $X^{(x)}$ for the PSSMP starting from $x\geq 0$.
Well-known examples of this kind of processes are: Bessel processes,
stable subordinators and  stable
L\'evy processes conditioned to stay positive.\\
In this paper, we are interested in the class of processes which
drift towards $+\infty$. Let $X^{(x)}$ be a PSSMP starting at $x\geq
0$ which drifts towards $+\infty$. We introduce the so-called future
infimum of $X^{(x)}$, by
\[
J^{(x)}_{t}\eqdef\inf_{s\geq t}X^{(x)}_{s},  \qquad\text{ for }
\quad t\geq 0.
\]
Note that the  future infimum process $J^{(x)}=(J^{(x)}_{t}, t\geq
0)$, is an increasing self-similar process with the same scaling
coefficient as $X^{(x)}$. It is clear that when the PSSMP $X^{(x)}$
starts from $x=0$, the process $J^{(0)}$ starts also from  $0$. When
the PSSMP $X^{(x)}$ starts from $x>0$, the future infimum $J^{(x)}$
starts from the global infimum, that is from $\inf_{t\geq
0}X^{(x)}_{t}$.
In both cases, the future infimum process $J^{(x)}$ tends to $+\infty$ as $t$ increases.\\
We are interested in describing the upper envelope at $0$ and at
$+\infty$ of the future infimum process for a large class of PSSMP
through integral tests and laws of iterated logarithm. As we will
note later, the same integral tests will also allow us to describe
the upper envelope at $0$ and $+\infty$ of the PSSMP $X^{(0)}$ in
the increasing case.\\
Khoshnevisan et al.~\cite{Kh}  studied the asymptotic behaviour of
the future infimum of some stochastic processes, in particular the
case of a Bessel process with index $d>2$. They obtained the following integral test:\\
\textit{Let $\phi(t)=\sqrt{t}\psi(t)$ be nondecreasing in $t>0$
and assume that $\psi(t)$ diverges to $+\infty$ as $t$ goes to
$+\infty$. If}
\[
\int^{+\infty}\Big(\psi(t)\Big)^{d-2}\exp\left\{-\psi^2(t)/2\right\}\ud
t<+\infty,
\]
\textit{then,  for all } $\epsilon>0$
\[
\p_{0}\left(J_t>(1+\epsilon)\phi(t), \textrm{ i.o.}, \textrm{ as }
t\to +\infty\right)=0.
\]
They also obtained the following law of the iterated logarithm,
\[
\limsup_{t\to +\infty}\frac{J_t}{\sqrt{2t\log\log(t)}}=1, \qquad
\p_0-\textrm{ a.s.}
\]
In this paper, we will give integral tests for $J$ in the general case.\\
Our arguments are based in the following representation of
self-similar Markov processes well-known as Lamperti representation.
Lamperti studied in detail the PSSMP in \cite{la}. In his main
result, Lamperti proved that any PSSMP starting from a strictly
positive state is a time-change of the exponential of a L\'evy
process. More precisely, let $X^{(x)}$ be a self-similar Markov
process started from $x>0$  that fulfills the scaling property for
some $\alpha>0$, then there exists $\xi=(\xi_{t}, t\geq 0)$ a L\'evy
process possibly killed at an independent exponential time, such
that
\begin{equation}\label{lamp}
X^{(x)}_{t}=x\exp\Big\{\xi_{\tau(tx^{-\alpha})}\Big\},\qquad 0\leq
t\leq x^{\alpha}I(\xi),
\end{equation}
where
\begin{equation*}
\tau_{t}=\inf\Big\{s\geq 0: I_{s}(\xi)>t\Big\},\quad
I_{s}(\xi)=\int_{0}^{s}\exp\big\{\alpha\xi_{u}\big\}\ud
u\quad\textrm{and}\quad I(\xi)=\lim_{t\to+\infty}I_{t}(\xi).
\end{equation*}
Lamperti raised the question of whether one can make sense of
$X^{(x)}$ started from $0+$. This problem was first solved for the
increasing case by Bertoin and Caballero \cite{beC}, this is when
its associated L\'evy process is a subordinator. In the same work,
they also computed the entrance law at $0$ which will be written
below. Later, Bertoin and Yor \cite{BeY} studied and computed the
entrance law in a more general case, when $\xi$ satisfies the
following condition
\[
(\textrm{H})\qquad 0<\e(\xi_{1})\leq
\e(|\xi_{1}|)<\infty \quad\textrm{ and } \quad \xi_{1}
\quad\textrm{is not arithmetic}.
\]
It was proved in \cite{BeY} that under condition (H), the process
$X^{(x)}$ converges in the sense of finite dimensional
distributions. If we denote by $\p_0$ its limit law then the
entrance law can be expressed as follows: for every measurable
function $f:\R_{+}\to\R_{+}$ and every $t>0$,
\begin{equation}\label{entrlaw}
\e_{0}\big(f(X_{t})\big)=\frac1m
\e\Big(I(\hat{\xi})^{-1}f\big(tI(\hat{\xi})^{-1}\big)\Big)
\end{equation}
where $m=E(\xi_{1})$ and $\hat{\xi}=-\xi$. Recently  Caballero and
Chaumont \cite{ChC} gave necessary and sufficient conditions for the
weak convergence of $X^{(x)}$ on the Skorokhod's space. In the
mentioned study , they give a path construction of $X^{(0)}$.\\
In this paper, we suppose that the limiting process $X^{(0)}$
exists in the sense of weak convergence on the Skorokhod's space
and that it satisfies that $\lim_{t\to +\infty}
X^{(0)}_{t}=+\infty$ which is equivalent, according to \cite{ChC}
to the fact that $\xi$ satisfies condition (H). Let $\mathcal{D}$
denote the  Skorokhod's space of c\`adl\`ag paths with real values
and defined on $[0, +\infty)$ and $\p$  a probability measure
defined on $\mathcal{D}$. We denote by $\p_{x}$, for $x> 0$ the
law, under $\p$, of the process $X^{(x)}$ defined above in
(\ref{lamp}) and by $\p_{0}$ the law, under $\p$, of the limiting
process $X^{(0)}$ whose entrance law is given by (\ref{entrlaw}).
With this notation we have that $(X, \p_{x})=(X^{(x)}, \p)$ for
$x\geq 0$. Throughout this work we will suppose that $\xi$ is a
L\'evy process satisfying condition (H).\\
Note that from the scaling property, the process $(X^{\alpha},
\p_{x})$, for $x\geq 0$ is a PSSMP whose scaling coefficient is
equal to $1$. Henceforth, without loss of generality we can assume
that
$\alpha=1$.\\
Some general results for the {\it lower} envelope of $X^{(x)}$ have
been established by Chaumont and Pardo in \cite{ChP}. These results
are based on the study of the last passage times of the process
$X^{(x)}$. Since its future infimum process $J^{(x)}$ can be seen as
the right inverse of the last passage times of $X^{(x)}$, it is not
difficult to deduce that we can replace $X^{(x)}$ by its future
infimum in all their results. In other words, we will obtain the
same integral tests for the {\it lower} envelope of $J^{(x)}$ at 0
$( \textrm{ when } x=0)$ and at $+\infty$ $(\textrm{ for all }x\geq
0)$. Motivated by the above and based on the study of the last
passage times for the PSSMP by Chaumont and Pardo, we will describe
the {\it upper} envelope of its
future infimum process.\\
For $x>0$, we consider $\hat{X}^{(x)}$, the dual process of the
PSSMP $X$ with respect to the Lebesgue measure. From Bertoin and Yor
\cite{BeY}, we know that $\hat{X}^{(x)}$ has a Lamperti
representation and is given by
\[
\hat{X}^{(x)}=\left(x\exp\Big\{\hat{\xi}_{\hat{\tau}(t/x)}\Big\},
0\leq t\leq xI\big(\hat{\xi}\big) \right),
\]
where
\[
\hat{\tau}_{t}=\inf\Big\{s\geq 0:I_{s}(\hat{\xi})>t\Big\},\qquad\
I_{s}(\hat{\xi})=\int_{0}^{s}\exp\Big\{\hat{\xi}_{t}\Big\}\ud
t\quad\textrm{and}\quad I(\hat{\xi})=\lim_{s\to
+\infty}I_{s}(\hat{\xi}).
\]
Note,  that  $xI(\hat{\xi})$ is the first time at which the
process $\hat{X}^{(x)}$ reaches the state $0$, that is
$xI(\hat{\xi})=\inf\{t: \hat{X}^{(x)}_{t}=0\}$.\\
For $y\geq 0$, we
define  $T_{y}=\inf\{t:\xi_{t}\geq y\}$, the first passage time of
the process $\xi$ over the state $y$,  and
$U(y)=\sup\{t:X^{(0)}_{t}\leq y\}$ the last passage time of the
processes $X^{(0)}$ below $y$. Since the processes $\xi$ and
$X^{(0)}$
drift towards $+\infty$, both random times are almost surely finite.\\
The following lemma proved by Chaumont and Pardo \cite{ChP}, gives a
path decomposition of the process $X^{(0)}$ reversed at time $U(x)$,
for $x>0$, and also determines the law of $U(x)$. This path
decomposition consists in splitting the path of $(X^{(0)}_{t},0\leq
t\leq U(x))$ at its last passage times.
\begin{lemma}Let $\Gamma=X^{(0)}_{U(x)^{-}}$. Then, the process time-reversed at
its last passage time below $x>0$,
$\hat{X}\eqdef(X^{(0)}_{(U(x)-t)^{-}}, 0\leq t\leq U(x))$ can be
described as
\[
\hat{X}=\Big(\Gamma\exp\Big\{\hat{\xi}_{\hat{\tau}(t/\Gamma)}\Big\},
0\leq t\leq U(x)\Big),
\]
where $\Gamma$ and $\hat{\xi}$ are independent. Moreover, let
$(x_{n})$ be a decreasing sequence which tends to $0$ and such that
$x_{1}=x$.  If we define $\hat{S}_y=\inf\{t: \hat{X}_{t}\leq y\}$,
for $y>0$, then we can describe the process $\hat{X}$ between the
passage times $\hat{S}_{x_{n}}$ and $\hat{S}_{x_{n+1}}$ as follows:
\[
\Big(\hat{X}_{\hat{S}_{x_{n}}+t}, 0\leq t\leq
\hat{S}_{x_{n+1}}-\hat{S}_{x_{n}}\Big)=\Big(\Gamma_{n}\exp\Big\{\hat{\xi}^{(n)}_{\hat{\tau}^{(n)}(t/\Gamma_{n})}\Big\},
0\leq t\leq H_{n}\Big), \quad n\geq 1,
\]
where the processes $\hat{\xi}^{(n)}$, $n\geq 1$ are independent
between themselves and have the same law as $\hat{\xi}$ and
\begin{align}
\hat{\tau}^{(n)}_{t}&= \inf\Big\{s: I_{s}\big(\hat{\xi}^{(n)}\big)>t\Big\},\\
H_{n}&
=\Gamma_{n}\int_{0}^{\hat{T}^{(n)}(\log(x_{n+1}/\Gamma_{n}))}\exp\Big\{\hat{\xi}^{(n)}_{s}\Big\}\ud
s,\\
\Gamma_{n+1}&=
\Gamma_{n}\exp\Big\{\hat{\xi}^{(n)}_{\hat{T}^{(n)}(\log(x_{n+1}/\Gamma_{n}))}\Big\},\quad
n\geq 1,\quad \Gamma_{1}=\Gamma,\\
\hat{T}^{(n)}_{z}&=\inf\Big\{ t:\hat{\xi}^{(n)}_{t}\leq z\Big\}.
\end{align}
For each $n\geq 1$, $\Gamma_{n}$ is independent of $\xi^{(n)}$ and
\begin{equation}\label{igusalt}
x^{-1}_{n}\Gamma_{n}\ed x^{-1}\Gamma,
\end{equation}
where the latter equality means that both variables have the same
distribution.
\end{lemma}
As a consequence, we have that for all $n\geq 1$,
\begin{equation}\label{suma}
U(x_{n})=\sum_{k\geq
n}\Gamma_{k}\int_{0}^{\hat{T}^{(k)}(\log(x_{k+1}/\Gamma_{k}))}\exp\Big\{\hat{\xi}^{(k)}_{s}\Big\}\ud
s,\quad \textrm{a.s.}
\end{equation}
On the other hand, since
\begin{equation}\label{iguallem}
\begin{split}
\Big(X^{(0)}_{(U(x_{n})-t)^{-}}, 0\leq t\leq U(x_{n})
\Big)&=\Big(\hat{X}_{\hat{S}_{x_{n}}+t}, 0\leq t\leq
U(x_{1})-\hat{S}(x_{n})\Big)\\
&=\Big(\Gamma_{n}\exp\Big\{\bar{\xi}^{(n)}_{\bar{\tau}^{(n)}(t/\Gamma_{n})}\Big\},
0\leq t\leq \Gamma_{n}I(\bar{\xi}^{(n)})\Big),
\end{split}
\end{equation}
where $\bar{\xi}^{(n)}$ has the same law as $\hat{\xi}$ and
$\bar{\tau}^{(n)}$ is the inverse of the exponential functional
$I_{s}(\bar{\xi}^{(n)})$. Then we also have that almost surely
\begin{equation}\label{igual}
U(x_{n})=\Gamma_{n} I\big(\bar{\xi}^{(n)}\big)\leq
x_{n}I\big(\bar{\xi}^{(n)}\big),\quad\textrm{where}\quad
I\big(\bar{\xi}^{(n)}\big)=\int_{0}^{\infty}\exp{\big\{\bar{\xi}^{(n)}_{t}\big\}}\ud
t,
\end{equation}
and, $\Gamma_{n}$ and $\bar{\xi}^{(n)}$ are independent.\\
Note that the process $(\hat{\xi}^{(n)}_{t}, 0\leq t\leq
\hat{T}^{(n)}(\log(x_{n+1}/\Gamma_{n}))$ is the same as the process
$\bar{\xi}^{(n)}$ killed  at
$\bar{T}^{(n)}(\log(x_{n+1}/\Gamma_{n}))$, where
$\bar{T}^{(n)}(x)=\inf\{ t:\bar{\xi}^{(n)}_{t}\leq x \}$ for $x\leq
0$. In fact, the process $\bar{\xi}^{(n)}$ can be described as
follows,
\begin{equation}\label{conca}
\bar{\xi}^{(n)}_{t}=\left\{\begin{array}{ll}
\hat{\xi}_{t}^{(n)} & \textrm{ if } t\in[0,\Sigma^{(n)}_{1}[,\\
\hat{\xi}_{t-\Sigma^{(n)}_{1}}^{(n+1)} & \textrm{ if } t\in[\Sigma^{(n)}_{1},\Sigma^{(n)}_{2}[,\\ \vdots & \\
\hat{\xi}_{t-\Sigma^{(n)}_{k}}^{(n+k)} & \textrm{ if } t\in[\Sigma^{(n)}_{k},\Sigma^{(n)}_{k+1}[,\\
\vdots &
\end{array} \right.
\end{equation}
where $\Sigma^{(n)}_{k}=\sum_{j=n}^{n+k-1}\hat{T}^{(j)}$ and
$\hat{T}^{(j)}=\hat{T}^{(j)}(\log(x_{j+1}/\Gamma_{j}))$.\\
Chaumont and Pardo proved in the same work that we have the same
properties for $x$ large (see Corollary 3 in \cite{ChP}). This will
be very useful to establish
our asymptotic results at $+ \infty$.\\
It is important to note that the law of $x^{-1}_{1}\Gamma$ is
related with the upward ladder height process $\sigma=(\sigma_{t},
t\geq 0)$ associated to $\xi$ (see Bertoin ~\cite{be} for a proper
definition). In fact, its law (see for instance Lemma 1 in
~\cite{ChP}) is the same as that of $\exp\{-UZ\}$, where $U$ and
$Z$ are independent random variables, $U$ is uniformly distributed
over $[0,1]$ and the law of $Z$ is given by
\[
\p(Z>u)=\e(\sigma_{1})^{-1}\int_{(u,\infty)}s\mu(\ud s), \qquad
u\geq 0,
\]
where $\mu$ is the L\'evy measure of $\sigma$. In particular, we can deduce that for all $y< x_{1}$, $\p(\Gamma>y)>0$.\\
The following result gives us integral tests at $0$ for the upper
envelope of $J^{(0)}$. This theorem means in particular that the
asymptotic behaviour of $J^{(0)}$ only depends on the tail behaviour
of the law of $\nu I(\hat{\xi})$ and this of $I(\hat{\xi})$, where
$\nu$ is independent of $I(\hat{\xi})$ and has the same distribution
as $x^{-1}_{1}\Gamma$. Note that the support of the law of $\nu$ is the interval $[0,1]$.\\
Let us define
\[
\bar{F}_{\nu}(t)\eqdef\p\Big(\nu I\big(\hat{\xi}\big)<
t\Big)\quad\textrm{and}\quad
\bar{F}(t)\eqdef\p\Big(I\big(\hat{\xi}\big)< t\Big)
\]
and denote by $\mathcal{H}_{0}$ the totality of positive
increasing functions $h(t)$ on $(0, \infty)$ that satisfy
\begin{itemize}
\item[i)] $h(0)=0$, and
\item[ii)] there exists $\beta\in(0,1)$ such that $\displaystyle\sup_{t<\beta}\displaystyle\frac{t}{h(t)}<\infty.$
\end{itemize}
\begin{theorem}
Let $h\in\mathcal{H}_{0}$.
\begin{itemize}
\item[i)] If
\[
\int_{0^{+}}\bar{F}_{\nu}\left(\frac{t}{h(t)}\right)\frac{\ud
t}{t}<\infty,
\]
then for all $\epsilon > 0$
\[
\p_{0}\Big(J_{t}>(1+\epsilon)h(t),\textrm{ i.o., as }t\to 0\Big)=0.
\]
\item[ii)]  If
\[
\int_{0^{+}}\bar{F}\left(\frac{t}{h(t)}\right)\frac{\ud
t}{t}=\infty,
\]
then for all $\epsilon > 0$
\[
\p_{0}\Big(J_{t}>(1-\epsilon)h(t),\textrm{ i.o., as }t\to 0\Big)=1.
\]
\end{itemize}
\end{theorem}
For the integral tests at $+\infty$, we define
$\mathcal{H}_{\infty}$, the totality of positive increasing
functions $h(t)$ on $(0, \infty)$ that satisfy
\begin{itemize}
\item[i)] $\lim_{t\to\infty}h(t)=\infty$, and
\item[ii)] there exists $\beta>1$ such that $\displaystyle\sup_{t>\beta}\displaystyle\frac{t}{h(t)}<\infty.$
\end{itemize}
Then the upper envelope of $J^{(x)}$ at $+\infty$ is given by the
following result.
\begin{theorem} Let $h\in\mathcal{H}_{\infty}$.
\begin{itemize}
\item[i)] If
\[
\int^{+\infty}\bar{F}_{\nu}\left(\frac{t}{h(t)}\right)\frac{\ud
t}{t}<\infty,
\]
then for all $\epsilon > 0$ and for all $x\geq 0,$
\[
\p_{x}\Big(J_{t}>(1+\epsilon)h(t),\textrm{ i.o., as }t\to
+\infty\Big)=0.
\]
\item[ii)]  If
\[
\int^{+\infty}\bar{F}\left(\frac{t}{h(t)}\right)\frac{\ud
t}{t}=\infty,
\]
then for all $\epsilon > 0$ and for all $x\geq 0$
\[
\p_{x}\Big(J_{t}>(1-\epsilon)h(t),\textrm{ i.o., as }t\to
+\infty\Big)=1.
\]
\end{itemize}
\end{theorem}
The rest of this paper is organized as follows. In section 2, we
state two Theorems that describe the lower envelope of the last
passage time process $U$ at 0 and at $+\infty$, respectively. In
section 3 we prove Theorems 1 and 2. Sections 4 and 5 are devoted to
the ``\textit{regular}" case and ``\textit{logregular}" case
respectively and some examples. In section
 6, we study the case of transient Bessel processes and finally in
 section 7 we  discuss the upper envelope of the increasing
 self-similar Markov processes.
\section{The lower envelope of the last passage times}
Let us recall the definition of the last passage time of $X^{(0)}$,
\[
U(x)=\sup\Big\{t\geq 0: \mbox{}X^{(0)}_{t}\leq
x\Big\}\qquad\textrm{ for } x\geq 0.
\]
From this definition,  we see that $U=(U(x), x\geq 0)$ is also an
increasing self-similar process whose scaling coefficient is the
inverse of the scaling coefficient of $X^{(0)}$. Since the process
$X^{(0)}$ starts at $0$ and drifts towards $+\infty$, we deduce that
the process $U$ also starts at $0$
and tends to infinity as $x$ increases.\\
In this section, we are interested in the study of the behaviour
of process $U$ at $0$ and at $+\infty$. As we will see in the
following section, the asymptotic behaviour of process $U$ is
related to the asymptotic behaviour of the future infimum of
$X^{(0)}$. In fact, we will see that the lower envelope of the
last passage time process
$U$ determines the upper envelope of  $J^{(0)}.$\\
The following result will give us integral tests at $0$
for the lower envelope of $U$.\\
Let us denote by $\mathcal{H}^{-1}_{0}$ the totality of positive
increasing functions $h(x)$ on $(0, \infty)$ that satisfy
\begin{itemize}
\item[i)] $h(0)=0$, and
\item[ii)] there exists $\beta\in(0,1)$ such that $\sup_{x<\beta}x^{-1}h(x)<\infty.$
\end{itemize}
\begin{theorem}
Let $h\in\mathcal{H}^{-1}_{0}$.
\begin{itemize}
\item[i)] If
\[
\int_{0^{+}}\bar{F}_{\nu}\left(\frac{h(x)}{x}\right)\frac{\ud
x}{x}<\infty,
\]
then for all $\epsilon > 0$
\[
\p\Big(U(x)<(1-\epsilon)h(x),\textrm{ i.o., as }x\to 0\Big)=0.
\]
\item[ii)] If
\[
\int_{0^{+}}\bar{F}\left(\frac{h(x)}{x}\right)\frac{\ud
x}{x}=\infty,
\]
then for all $\epsilon > 0$
\[
\p\Big(U(x)<(1+\epsilon)h(x),\textrm{ i.o., as }x\to 0\Big)=1.
\]
\end{itemize}
\end{theorem}
\textit{Proof: } We first prove the convergent part. Let $(x_{n})$
be a decreasing sequence of positive numbers which converges to
$0$ and let us define the
events $A_{n}=\{U(x_{n+1})<h(x_{n})\}$.\\
Now, we choose $x_{n}=r^{n}$, for $r<1$. From the first Borel
Cantelli's Lemma, if $\sum_{n}\p(A_{n})<\infty$, it follows
\[
U\big(r^{n+1}\big)\geq h\big(r^{n}\big)\qquad \p-\textrm{ a.s.,}
\]
for all large $n$. Since the function $h$ and the process $U$ are
increasing, we have
\[
U(x)\geq h(x)\qquad \textrm{for}\quad r^{n+1}\leq x\leq r^{n}.
\]
From (\ref{igual}), we get the following inequality
\[
\begin{split}
\sum_{n}\p\Big(U\big(r^{n}\big)<h\big(r^{n+1}\big)\Big)&\leq
\int_{1}^{\infty}\p\Big(r^{t}\nu I\big(\hat{\xi}\big)<h\big(r^{t}\big)\Big)\ud t\\
&=-\frac{1}{\log{r}}\int^{r}_{0}\bar{F}_{\nu}\left(\frac{h(x)}{x}\right)\frac{\ud
x}{x}.
\end{split}
\]
From our hypothesis, this last integral is finite. Then from the
above discussion, there exist $x_{0}$ such that for every $x\geq
x_{0}$
\[
U(x)\geq r^{2}h(x),\qquad\textrm{for all }\quad r<1.
\]
Clearly, this implies that
\[
\p\Big( U(x)<r^{2}h(x), \textrm{ i.o., as }x\to 0\Big)=0,
\]
which proves part $(i)$.\\
Now we prove the divergent part. First, note that when $\xi$ has
no positive jumps the process $U$ is like the ones considered in
Watanabe's work, that is   $U$ is an increasing self-similar
process with independent increments, but in the general case the
process $U$ does not have this property. The decomposition
(\ref{suma}) and the a.s. equality in (\ref{igual}) will allow us
to
extend the arguments used by Watanabe to our case.\\
Now, we assume that $h$ satisfies
\[
\int_{0^{+}}\bar{F}\left(\frac{h(x)}{x}\right)\frac{\ud
x}{x}=\infty.
\]
Let us take, again $x_{n}=r^n$ for $r<1$, and define the events
\[
 \quad C_{n}=\Big\{U(x)<r^{-2}h(x), \textrm{ for some
}x\in (0, r^{n})\Big\}.
\]
Note that the family $(C_{n})$ is decreasing, then
\[
C=\bigcap_{n\geq1} C_{n}=\Big\{U(x)<r^{-2}h(x), \textrm{ i.o., as
}x\to 0\Big\}.
\]
If we prove that $\lim\p(C_{n})>0$, then since $X^{(0)}$ is a Feller
process and by Blumenthal's 0-1 law we will have that
\[
\p\Big( U(x)<r^{-2}h(x),\textrm{ i.o., as }x\to 0\Big)=1,
\]
which will prove part $(ii)$.\\
In this direction, we define the following events.  For $n\leq m-1$,
\[
D_{(n , m)}=\Big\{r^{j+1}\bar{I}_{(j+1,m+1)}\geq h(r^j),\textrm{ for
all } n\leq j\leq m-1\Big\},\] and for $r<k<1$ and $n\leq m-2$
\begin{align*}
E_{(n,m-1)}&=\Big\{r^{j+1}\bar{I}_{(j+1,m)}+ r^{j+1}R_{(j+1,
m)}\bar{I}_{(m, m+1)}\geq h(r^j),\textrm{ for all } n\leq j\leq
m-2\Big\}\quad\textrm{and}\\
E^{(k)}_{(n,m-1)}&=\Big\{r^{j+1}\bar{I}_{(j+1,m)}+ r^{j+1}R_{(j+1,
m)}\bar{I}^{(k)}_{m}\geq h(r^j),\textrm{ for all } n\leq j\leq
m-2\Big\},
\end{align*}
where
\begin{align*}
\bar{I}_{(j+1,m+1)}&=\int_{0}^{\bar{T}^{(j+1)}(\log(r^{m+1}/\Gamma_{j+1}))}\exp
\Big\{\bar{\xi}_{s}^{(j+1)}\Big\}\ud s ,\\
\bar{I}^{(k)}_{m}&=\int_{0}^{\bar{T}^{(m)}(\log(r^{m+1}/kr^{m}))}\exp
\Big\{\bar{\xi}_{s}^{(m)}\Big\}\ud s \qquad \textrm{ and}\\
R_{(j+1,m)}&=\exp\Big\{\bar{\xi}^{(j+1)}_{\bar{T}^{(j+1)}(\log(r^{m}/\Gamma_{j+1}))}\Big\},
\end{align*}
and for $n\leq j\leq m-1$, $\bar{\xi}^{(j+1)}$ is a L\'evy process defined as in (\ref{conca}).\\
From the definition of $\bar{\xi}^{(j+1)}$, we can deduce that for
$j<m$
\begin{align*}
\bar{\xi}^{(m)}&=\Big(\bar{\xi}^{(j+1)}_{\bar{T}^{(j+1)}(\log(r^{m}/\Gamma_{j+1}))+t}-\bar{\xi}^{(j+1)}_{\bar{T}^{(j+1)}(\log(r^{m}/\Gamma_{j+1}))},
t\geq 0\Big)\quad\textrm{ and }\\
\Gamma_{m}&=\Gamma_{j+1}\exp\Big\{\bar{\xi}^{(j+1)}_{\bar{T}^{(j+1)}(\log(r^{m}/\Gamma_{j+1}))}\Big\},
\end{align*}
then it is straightforward that
\[
\bar{T}^{(j+1)}\big(\log(r^{m+1}/\Gamma_{j+1})\big)=\bar{T}^{(j+1)}\big(\log(r^{m}/\Gamma_{j+1})\big)+\inf\Big\{t\geq
0;\mbox{} \bar{\xi}^{(m)}_{t}\leq \log(r^{m+1}/\Gamma_{m})\Big\}.
\]
The above decomposition allows us to determine the following
identity
\begin{equation}\label{iguaint}
\bar{I}_{(j+1,m+1)}=\bar{I}_{(j+1,m)}+ R_{(j+1, m)}\bar{I}_{(m,
m+1)}.
\end{equation}
In the same way we can also get that,
\begin{equation}\label{iguaint2}
I(\bar{\xi}^{(j+1)})=\bar{I}_{(j+1,m+1)}+R_{(j+1,
m+1)}I(\bar{\xi}^{(m+1)}).
\end{equation}
By Lemma 1 and the decomposition (\ref{conca}), it follows that
$I\big(\bar{\xi}^{(m+1)}\big)$ is independent of
$(\bar{I}_{(j+1,m+1)},R_{(j+1, m+1)})$ and distributed as
$I\big(\hat{\xi}\big)$.\\
From (\ref{iguaint}) and since
\[
\Big\{r^{m}\bar{I}_{(m,m+1)}\geq
h(r^{m-1})\Big\}\subset\big\{\Gamma_{m}> r^{m+1}\big\},
\] we conclude that
\[
D_{(n,m)}=E_{(n,m-1)}\bigcap \big\{r^{m}\bar{I}_{(m,m+1)}\geq
h(r^{m-1})\big\}\bigcap\big\{\Gamma_{m}> r^{m+1}\big\}.
\]
Now, for $n\leq m-1$, we define
\[
H(n,m)=\p\Big(E^{(k)}_{(n,m-1)},r^{m}\bar{I}^{(k)}_{m}\geq
h(r^{m-1}),\Gamma_{m}> r^{m}k \Big).
\]
On the event $\{\Gamma_{m}>r^mk\}$, we have that
$\bar{I}^{(k)}_{m}\leq \bar{I}_{(m,m+1)}$. Hence since $k>r$, we
 deduce that $\p(D_{(n, m)})\geq
H(n,m)$.\\
For our purpose, we will prove that there exist ($n_{l}$) and
($m_{l}$), two increasing sequences such that $0\leq n_{l}\leq
m_{l}-1$, and $n_{l}, m_{l}$ go to $ \infty$ and $H(n_{l}, m_{l})$
tends to $0$ as $l$ goes to infinity. In this direction, we define
the events
\[
B_{n}=\Big\{r^{n+1}I\big(\bar{\xi}^{(n+1)}\big)<h(r^{n})\Big\}.
\]
If we suppose the contrary, this is that there
exists $\delta>0$ such that $H(n, m)\geq \delta$ for all
sufficiently large integers $m$ and $n$, we see from identity
(\ref{iguaint2}) that
\[
\begin{split}
1&\geq \p\left(\bigcup_{m=n+1}^{\infty}B_{m}\right)\geq
\sum_{m=n+1}^{\infty}\p\left(B_{m}\bigcap\left(\bigcap_{j=n}^{m-1}B^c_{j}\right)\right)\\
&=\sum_{m=n+1}^{\infty}\p\left(r^{m+1}I\big(\bar{\xi}^{(m+1)}\big)<h(r^m),\bigcap_{j=n}^{m-1}\Big\{r^{j+1}I(\bar{\xi}^{(j+1)})\geq
h(r^j)\Big\}\right)\\
&\geq\sum_{m=n+1}^{\infty}\p\left(r^{m+1}I\big(\bar{\xi}^{(m+1)}\big)<h(r^m)\right)\p\big(D_{(n,m)}\big)\\
&\geq\sum_{m=n+1}^{\infty}\p\Big(r^{m+1}I\big(\bar{\xi}^{(m+1)}\big)<h(r^m)\Big)H(n,m)\geq
\delta\sum_{m=n+1}^{\infty}\p\Big(r^{m+1}I\big(\hat{\xi}\big)<h(r^m)\Big),
\end{split}
\]
but this last sum diverges, since
\[
\begin{split}
\sum_{m=n+1}^{\infty}\p\big(r^{m+1}I\big(\hat{\xi}\big)<h(r^m)\big)&\geq
\int_{n+1}^{\infty}\p\Big(r^tI\big(\hat{\xi}\big)<h(r^t)\Big)\ud
t\\
&=-\frac{1}{\log
r}\int_{0}^{r^{n+1}}\bar{F}\left(\frac{h(x)}{x}\right)\frac{\ud
x}{x}.
\end{split}
\]
Hence our assertion is true.\\
Next, we denote $\p(I(\hat{\xi})\in \ud x)=\mu(\ud x)$ and
$\p(I_{r/k}\in \ud x)=\bar{\mu}(\ud x)$ for $k>r$, where
$I_{r/k}=\int_{0}^{\hat{T}_{\log(r/k)}}\exp\{\hat{\xi}_{s}\} \ud s$,
and we define
\[
\rho_{n_{l},m_{l}}(x)=\p\left(\bigcap
_{j=n_{l}}^{m_{l}-2}\Big\{r^{j+1}\bar{I}_{(j+1,m_{l})}+r^{j+1}R_{j+1,m_{l}}x\geq
r^{-1}h(r^j)\Big\}, \Gamma_{m_{l}}>kr^{m_{l}}\right),
\]
and
\[
G(n_{l},m_{l})=\p\left(
\bigcap_{j=n_l}^{m_l-1}\Big\{r^{j+1}I(\bar{\xi}^{(j+1)})\geq
h(r^j)\Big\}, \Gamma_{m_{l}}>kr^{m_{l}}\right).
\]
Note that $\rho_{n_l,m_l}(x)$ is increasing in $x$.\\
Hence, $H(n_l,m_l)$ and $G(n_l,m_l)$ are expressed as follows
\begin{align*}
H(n_l,m_l)&=\int_{r^{-m_l}h(r^{m_l-1})}^{\infty}\bar{\mu}(\ud
x)\rho_{n_l,m_l}(x)\quad\textrm{and}\\
G(n_l,m_l)&=\int_{r^{-m_l}h(r^{m_l-1})}^{\infty}\mu(\ud
x)\rho_{n_l,m_l}(x).
\end{align*}
The equality for $H(n_l, m_l)$ is evident since the random
variable $\bar{I}^{(k)}_{m}$ is independent from
$\big\{\Gamma_{m_l},(\bar{I}_{(j+1, m_{l})}, R_{(j+1, m_l)};
n_{l}\leq j\leq m_{l}-2 )\big\}$. To show the second one, we use
(\ref{iguaint2}) in the following form
\[
I\big(\bar{\xi}^{(j+1)}\big)=\bar{I}_{(j+1,m_{l})}+R_{(j+1,m_{l})}I\big(\bar{\xi}^{(m_{l})}\big),
\]
and the independence between $I\big(\bar{\xi}^{(m_l)}\big)$ and
$\big\{\Gamma_{m_l},(\bar{I}_{(j+1, m_{l})}, R_{(j+1, m_l)};
n_{l}\leq j\leq m_{l}-2 )\big\}$.\\
In particular, it follows that for $l$ sufficiently large
\[
H(n_{l}, m_{l}) \geq
\rho_{n_{l},m_{l}}(N)\int_{N}^{\infty}\bar{\mu}(\ud x)\qquad
\textrm{ for } \quad N\geq rC,
\]
where $C=\sup_{x\leq \beta} x^{-1}h(x)$.\\
Since $H(n_{l}, m_{l})$ converges to $0$, as $l$ goes to $+\infty$
and $\bar{\mu}$ does not depend on $l$, then $\rho_{n_{l},
m_{l}}(N)$ also converges to
$0$ when $l$ goes to $+\infty$, for every $N\geq rC$.\\
On the other hand, we have
\[
G(n_{l}, m_{l})\leq \rho_{n_{l}, m_{l}}(N)\int_{0}^{N}\mu(\ud
 x)+\int_{N}^{\infty}\mu(\ud x),
\]
then, letting $l$ and  $N$ go to infinity, we get that $G(n_l,
m_l)$ goes to $0$.\\
Note that the set $C_{n_l}$ satisfies
\[
\p(C_{n_l})\geq 1-\p\left(r^{j+1}I(\bar{\xi}^{(j+1)})\geq h(r^j),
\textrm{ for all } n_l\leq j\leq m_l -1\right)
\]
and it is not difficult to see that
\[
\p\left(r^{j+1}I(\bar{\xi}^{(j+1)})\geq h(r^j), \textrm{ for all }
n_l\leq j\leq m_l -1\right)\leq \p\big(\Gamma_{m_l}\leq
kr^{m_l})+G(n_{l}, m_{l}).
\]
Then,
\[
\p(C_{n_l})\geq\p\big(\Gamma_{m_l}> kr^{m_l})-G(n_{l}, m_{l}),
\]
and since $\p\big(\Gamma_{m_l}> kr^{m_l})=\p\big(\Gamma> kr)>0$ (see
Lemma 1 and the properties of $\Gamma$ in Section 1), we conclude
that $\lim\p(C_n)>0$ and with this we finish the proof.\QED For the
integral test at $+\infty$, we define $\mathcal{H}^{-1}_{\infty}$
the totality of positive increasing functions $h(x)$ on $(0,
\infty)$ that satisfy
\begin{itemize}
\item[i)] $\lim_{x\to+\infty}h(x)=+\infty$, and
\item[ii)] there exists $\beta>1$ such that $\sup_{x>\beta}x^{-1}h(x)<\infty.$
\end{itemize}
\begin{theorem}
Let $h\in\mathcal{H}^{-1}_{\infty}$.
\begin{itemize}
\item[i)] If
\[
\int^{+\infty}\bar{F}_{\nu}\left(\frac{h(x)}{x}\right)\frac{\ud
x}{x}<\infty,
\]
then for all $\epsilon > 0$
\[
\p\Big(U(x)<(1-\epsilon)h(x),\textrm{ i.o., as }x\to +\infty\Big)=0.
\]
\item[ii)] If
\[
\int^{+\infty}\bar{F}\left(\frac{h(x)}{x}\right)\frac{\ud
x}{x}=\infty,
\]
then for all $\epsilon > 0$
\[
\p\Big(U(x)<(1+\epsilon)h(x),\textrm{ i.o., as }x\to +\infty\Big)=1.
\]
\end{itemize}
\end{theorem}
\textit{Proof:} The proof is very similar to that in Theorem 3.
First, note that we have the same results as Lemma 1 for $x$ large
(see Corollary 3 of Chaumont and Pardo \cite{ChP}), then we get the
integral test following the same arguments for the proof of $(i)$
and $(ii)$ for the sequence $x_{n}=r^{n}$, for $r>1$, and noticing
that if we define
\[
\begin{split}
C_{n}&=\Big\{U(x)<h_{r}(x), \textrm{ for some } x\in (r^{n},
+\infty)\Big\}\\
&=\Big\{J^{(0)}_{t}>h_{r}^{-1}(t), \textrm{ for some }
t\in(U(r^n), +\infty)\Big\},
\end{split}
\]
where $h_{r}(t)=r^{2}h(t)$, then the event $C=\cap_{n\geq 1}
C_{n}$ is in the upper-tail sigma-field
$\cap_{t}\sigma\{X_{s}^{(0)}\mbox{}: \mbox{}s\geq t\}$ which is
trivial.\QED
\section{Proof of Theorems 1 and 2}
\textit{Proof of Theorem 1:} Let $(x_{n})$ be a decreasing sequence
which converges to $0$. We define the events
$A_{n}=\big\{\textrm{There exists }t\in [U(x_{n+1}),
U(x_{n})]\textrm{ such that } J^{(0)}_{t}>h(t)\big\}$. From the fact
that $U(x_{n})$ tends to $0$, a.s. when $n$ goes to $+\infty$, we
see
\[
\Big\{J^{(0)}_{t}>h(t), \textrm{ i.o., as }t\to
0\Big\}=\limsup_{n\to +\infty} A_{n}.
\]
Since $h$ is an increasing function and $J^{(0)}_{U(x_{n})}\geq
x_{n}$ a.s., the following inclusions hold
\begin{equation}\label{inclusion}
\Big\{x_{n}>h\big(U(x_{n})\big) \Big\}\subset A_{n}\subset
\Big\{x_{n}>h\big(U(x_{n+1})\big)\Big\}.
\end{equation}
Now, we prove the convergent part. We choose $x_{n}=r^{n}$, for
$r<1$ and $h_{r}(t)=r^{-2}h(t)$. Since $h$ is increasing, we deduce
that
\[
\sum_{n} \p\Big(r^n>h_{r}\big(U(r^{n+1})\big)\Big)\leq
-\frac{1}{\log r}\int_{0}^{r}\p\Big(t>h\big(U(t)\big)\Big)\frac{\ud
t}{t}.
\]
Replacing $h$ by $h_r$ in (\ref{inclusion}), we see that we can
obtain our result if
\[
 \int_{0}^{r}\p\Big(t>h\big(U(t)\big)\Big)\frac{\ud t}{t}<\infty.
\]
From elementary calculations, we deduce that
\[
\int_{0}^{r}\p\Big(t>h\big(U(t)\big)\Big)\frac{\ud t}{t}=\e \left(
\int_{0}^{h^{-1}(r)}\ind_{\big\{t/r<\nu
I(\hat{\xi})<t/h(t)\big\}}\frac{\ud t}{t}\right),
\]
where $h^{-1}(s)=\inf\{t>0, h(t)>s\}$, the right inverse function of
$h$. Then, this integral converges if
\[
\int_{0}^{h^{-1}(r)}\p\left(\nu
I\big(\hat{\xi}\big)<\frac{t}{h(t)}\right)\frac{\ud t}{t} <\infty.
\]
This proves part $(i)$.\\
Next, we prove the divergent case. We suppose that $h$ satisfies
\[
\int_{0^+}\bar{F}\left(\frac{t}{h(t)}\right)\frac{\ud
t}{t}=\infty.
\]
Take, again, $x_{n}=r^{n}$, for $r<1$ and  note that,
\[
\begin{split}
B_{n}&=\bigcup_{m=n}^{\infty}A_{m}=\big\{\textrm{There exist }t\in
(0, U(r^n)]\textrm{ such that } J^{(0)}_{t}>h_{r}(t)\big\}\\
&=\big\{\textrm{There exist }x\in (0, r^n]\textrm{ such that }
U(x)<h^{-1}_{r}(x)\big\}
\end{split}
\]
where $h_{r}(t)=rh(t)$ and $h^{-1}_{r}$ its right inverse function.
Hence, by analogous arguments to the proof of Theorem 3 part $(ii)$
it is enough to prove that $\lim\p(B_{n})>0$ to obtain our result.
With this purpose, we will follow the proof of Theorem 3.\\
From inclusion (\ref{inclusion}) and the a.s. inequality in
(\ref{igual}) we see
\[
\p(B_{n})\geq 1-\p\Big(r^{j}\leq
rh\big(r^{j}I(\bar{\xi}^{(j)})\big), \textrm{ for all } n\leq j\leq
m-1 \Big),
\]
where $m$ is chosen arbitrarily $m\geq n+1$.\\
Now, we define the events
\[
C_{n}=\bigg\{r^n>rh\Big(r^{n}I\big(\bar{\xi}^{(n)}\big)\Big)\bigg\},
\]
and we will prove that $\sum \p(C_{n})=\infty.$  Since the function
$h$ is increasing, it is straightforward that
\[
\sum_{n\geq 1} \p(C_{n}) \geq
\int_{0}^{+\infty}\p\left(r^{t}>h\Big(r^t
I\big(\hat{\xi}\big)\Big)\right)\ud t=-\frac{1}{\log
r}\int_{0}^{1}\p\Big(t>h\big(tI\big(\hat{\xi}\big)\big)\Big)\frac{\ud
t}{t}.
\]
Hence, it is enough to prove that this last integral is infinite. In
this direction, we have that
\[
\int_{0}^{r}\p\Big(t>h\big(tI\big(\hat{\xi}\big)\big)\Big)\frac{\ud
t}{t} =\e \left( \int_{0}^{h^{-1}(r)}\ind_{\big\{t/r<
I(\hat{\xi})<t/h(t)\big\}}\frac{\ud t}{t}\right).
\]
On the other hand, we see that
\[
\begin{split}
\int_{0}^{h^{-1}(r)}\p\left(
I\big(\hat{\xi}\big)<\frac{t}{h(t)}\right)\frac{\ud t}{t}=
&\int_{0}^{h^{-1}(r)}\p\left(\frac{t}{r}<
I\big(\hat{\xi}\big)<\frac{t}{h(t)}\right)\frac{\ud
t}{t}\\&+\int_{0}^{h^{-1}(r)}\p\left(
I\big(\hat{\xi}\big)<\frac{t}{r}\right)\frac{\ud t}{t},
\end{split}
\]
and since $e^{-1}\hat{T}_{1}\leq I(\hat{\xi})$ a.s., then
\[
\begin{split}
\int_{0}^{h^{-1}(r)}\p\left(
I\big(\hat{\xi}\big)<\frac{t}{r}\right)\frac{\ud t}{t}&\leq
\e\left(\log^{+}
\frac{h^{-1}(r)}{r}I\big(\hat{\xi}\big)^{-1}\right)\\
&\leq 1+\log^{+}\frac{h^{-1}(r)}{r}+\e\left(|\log
\hat{T}_{1}|\right),
\end{split}
\]
which is clearly finite from our assumptions. Then, we deduce that
\[
\e \left( \int_{0}^{h^{-1}(r)}\ind_{\big\{t/r<
I(\hat{\xi})<t/h(t)\big\}}\frac{\ud t}{t}\right)=\infty,
\]
and hence $\sum\p(C_{n})=\infty$.\\
Next following the same notation as in the proof of Theorem 3, we
define the following events. For $n\leq m-1$
\[
D_{(n,m)}=\Big\{r^{j}\leq rh\big(r^{j}\bar{I}_{(j,m)}\big), \textrm{
for all }n\leq j\leq m-1\Big\},
\]
and, for $r<k<1$ and $n\leq m-2$
\begin{align*}
E_{n,m-1}&=\Big\{r^j\leq rh\big(r^j\bar{I}_{(j,
m-1)}+r^{j}R_{(j,m-1)}\bar{I}_{(m-1,m)}\big), \textrm{ for all
}n\leq j\leq m-2\Big\}\quad \textrm{and }\\
E^{(k)}_{n,m-1}&=\Big\{r^j\leq rh\big(r^j\bar{I}_{(j,
m-1)}+r^{j}R_{(j,m-1)}\bar{I}^{(k)}_{m-1}\big), \textrm{ for all
}n\leq j\leq m-2\Big\}.
\end{align*}
Again, we have that
\[
D_{(n,m)}=E_{(n,m-1)}\bigcap \Big\{r^{m-1}\leq
rh\big(r^{m-1}\bar{I}_{(m-1,m)}\big)\Big\}\bigcap\big\{\Gamma_{m-1}>r^{m}k\big\}.
\]
Now, for $n\leq m-1$, we define
\[
H(n,m)=\p\Big(E^{(k)}_{(m,m-1)}, r^{m-1}\leq
rh\big(r^{m-1}\bar{I}_{(m-1,m)}\big), \Gamma_{m-1}>r^{m}k\Big).
\]
Since $k>r$, we deduce that $\p(D_{(n,m)})\geq H(n,m)$. Then
similarly as in the proof of Theorem 3, we will prove that there
exist $(n_l)$ and $(m_{l})$, two increasing sequences such that
$0\leq n_{l}\leq m_{l}-1$, and $n_{l}, m_{l}$ go to $+\infty$ and
$H(n_l,m_l)$ tends to $0$ as $l$ goes to infinity.\\
We suppose the contrary, i.e., there exist $\delta>0$ such that
$H(n, m)\geq \delta$ for all sufficiently large integers $m$ and
$n$, hence
\[
\begin{split}
1&\geq \p\left(\bigcup_{m=n+1}^{\infty}C_{m}\right)\geq
\sum_{m=n+1}^{\infty}\p\left(C_{m}\bigcap\left(\bigcap_{j=n}^{m-1}C^{c}_{j}\right)\right)\\
&\geq\sum_{m=n+1}^{\infty}\p\Big(r^{m}>rh\Big(r^{m}I\big(\bar{\xi}^{(m)}\big)\Big)\Big)\p(D_{(n,m)})\\
&\geq\sum_{m=n+1}^{\infty}\p\Big(r^{m}>rh\Big(r^{m}I\big(\bar{\xi}^{(m)}\big)\Big)\Big)H(n,m)\geq
\delta
 \sum_{m=n+1}^{\infty}\p\big(C_{m}\big),\\
\end{split}
\]
but since $\sum \p(C_{n})$ diverges, we see that our assertion is
true.\\
Next, we define
\[
\rho_{n_{l},
m_{l}}(x)=\p\left(\bigcap_{j=n_{l}}^{m_{l}-2}\Big\{r^j\leq
rh\big(r^j\bar{I}_{(j, m-1)}+r^{j}R_{(j,m-1)}x\big)\Big\},
\Gamma_{m_{l}-1}>kr^{m_l-1}\right)
 \]
and
\[
G(n_{l}, m_{l})=\p\left(\bigcap_{j=n_{l}}^{m_{l}-1}\Big\{r^{j}\leq
rh\Big(r^{j}I\big(\bar{\xi}^{(j)}\big)\Big)\Big\},
\Gamma_{m_{l}-1}>kr^{m_l-1}\right).
\]
Since $h$ is increasing, we see that $\rho_{n_{l}, m_{l}}(x)$ is
increasing in $x$.\\
Again, we express $H(n_l,m_l)$ and $G(n_{l}, m_{l})$ as follows
\begin{align*}
H(n_l,m_l)&=\int_{0}^{+\infty}\bar{\mu}(\ud
x)\ind_{\big\{h(r^{m_{l}-1}x)\geq r^{m_{l}}\big\}}\rho_{n_{l},
m_{l}}(x)\quad \textrm{and},\\
G(n_l,m_l)&=\int_{0}^{+\infty}\mu(\ud
x)\ind_{\big\{h(r^{m_{l}-1}x)\geq r^{m_{l}}\big\}}\rho_{n_{l},
m_{l}}(x).
\end{align*}
In particular, we get that for $l$ sufficiently large
\[
H(n_l,m_l)\geq \rho_{n_{l}, m_{l}}(N)\int_{N}^{+\infty}\bar{\mu}(\ud
x)\rho_{n_{l}, m_{l}}(x)\qquad{for }\quad N\geq rC,
\]
where $C=\sup_{x\leq \beta}x/h(x)$. Hence following the same
arguments of the proof of Theorem 3, it is not difficult to see that
$G(n_l,m_l)$ goes to $0$ as $l$ goes to infinity and that
\[
\lim_{l\to+\infty}1-\p\Big(r^{j+1}\leq
rh\big(r^{j+1}I(\bar{\xi}^{(j+1)})\big), \textrm{ for all }
n_{l}\leq j\leq m_{l}-1 \Big)>0.
\]
Then, we conclude that  $\lim\p(B_{n})>0$ and with this we finish
the proof. \QED \noindent \textit{Proof of Theorem 2:} We first
consider the case where $x=0$. In this case the proof of the tests
at $+\infty$ is almost the same as that of the tests at $0$. It is
enough to apply the same arguments to the sequence $x_{n}=r^{n}$,
for $r>1$.\\
Now, we prove $(i)$ for any $x>0$. Let $h\in \mathcal{H}_{\infty}$
such that
$\int^{+\infty}\bar{F}_{\nu}\left(\frac{t}{h(t)}\right)\frac{\ud
t}{t}$ is finite. Let $x>0$ and $S_{x}=\inf\{t\geq
0\mbox{}:\mbox{}X^{(0)}_{t}\geq x\}$ and note by $\mu_{x}$ the law
of $X^{(0)}_{S_{x}}$. Since clearly
\[
\int^{+\infty}\bar{F}_{\nu}\left(\frac{t}{h(t-S_{x})}\right)\frac{\ud
t}{t}<\infty,
\]
from the Markov property at time $S_{x}$, we have
for all $\epsilon>0$
\begin{equation}\label{lliMark}
\begin{split}
&\p_{0}\Big(J_{t}>(1+\epsilon)h(t-S_{x}), \textrm{ i. o., as } t\to
\infty\Big)\\
=&\int_{[x,+\infty)}\p_{y}\Big(J_{t}>(1+\epsilon)h(t), \textrm{ i.
o., as } t\to \infty\Big)\mu_{x}(\ud y)=0.
\end{split}
\end{equation}
If $x$ is an atom of $\mu_{x}$, then equality (\ref{lliMark}) shows
that
\[
\p\Big(J^{(x)}_{t}>(1+\epsilon)h(t), \textrm{ i. o., as } t\to
\infty\Big)=0
\]
and the result is proved. Suppose that $x$ is not an atom of
$\mu_{x}$. Recall from Section 1, that $\log(x^{-1}_{1}\Gamma)$ is
the limit in law of the overshoot process
$\hat{\xi}_{\hat{T}_{x}}-x$, as $x\to+\infty$. So, it follows from
\cite{ChC}, Theorem 1 that
$X^{(0)}_{S_{x}}\ed\frac{xx_{1}}{\Gamma}$, and since $\p(\Gamma>z)$
for $z<x_{1}$, we have for any $\alpha>0$, $\mu_{x}(x, x+\alpha)>0$.
Hence (\ref{lliMark}) shows that there exists $y>x$ such that
\[
\p\Big(J^{(y)}_{t}>(1+\epsilon)h(t), \textrm{ i. o., as } t\to
\infty\Big)=0,
\]
for all $\epsilon>0$. The previous allows us to conclude part $(i)$.\\
Part $(ii)$ can be proved in the same way.\QED In some cases, it
will prove complicated to find sharp estimations of the tail
probability of $\nu I\big(\hat{\xi}\big)$, given that we will not
have enough information about the distribution of $\nu$. However, if
we can determine the law of $I\big(\hat{\xi}\big)$ then,  by
(\ref{entrlaw}) we will also determine the law of $X^{(0)}_{1}$ and
sometimes it will be possible to have sharp estimations of its tail
probability. For this reason, we will give another integral test for
the
convergence cases in Theorems 1 and 2, in terms of the tail probability of $X^{(0)}_{1}$. \\
Let us define
\[
G(t)=\p_{0}\big(X_{1}>t\big).
\]
\begin{corollary}
\begin{itemize}
\item[i)] Let $h\in \mathcal{H}_{0}$. If
\[
\int_{0^{+}}G\left(\frac{h(t)}{t}\right)\frac{\ud t}{t}<\infty,
\]
then for all $\epsilon > 0$
\[
\p\Big(J^{(0)}_t>(1+\epsilon)h(t),\textrm{ i.o., as }t\to 0\Big)=0.
\]
\item[ii)]Let $h\in \mathcal{H}_{\infty}$. If
\[
\int^{+\infty}G\left(\frac{h(t)}{t}\right)\frac{\ud t}{t}<\infty,
\]
then  and for all $\epsilon > 0$
\[
\p\Big(J^{(x)}_t>(1+\epsilon)h(t),\textrm{ i.o., as }t\to
\infty\Big)=0.
\]
\end{itemize}
\end{corollary}
\textit{Proof:} The proof of this Corollary is consequence of the
following inequality. By the scaling property,
\[
\bar{F}_{\nu}\big(t/h(t)\big)=\p\big(U(1)<t/h(t)\big)=\p\big(J^{(0)}_{1}>h(t)/t\big)\leq
\p_{0}\big(X_{1}>h(t)/t\big),
\]
and then applying Theorem 1 part (i), we obtain the desired
result.\QED In the same way, we can obtain another integral test for
the convergence cases in Theorem 3 and 4, in terms of $G$.
\begin{corollary}
\begin{itemize}
\item[i)] Let $h\in\mathcal{H}^{-1}_{0}$. If
\[
\int_{0^{+}}G\left(\frac{x}{h(x)}\right)\frac{\ud x}{x}<\infty,
\]
then for all $\epsilon > 0$
\[
\p\Big(U(x)<(1-\epsilon)h(x),\textrm{ i.o., as }x\to 0\Big)=0.
\]
\item[ii)]Let $h\in\mathcal{H}^{-1}_{\infty}$. If
\[
\int^{+\infty}G\left(\frac{x}{h(x)}\right)\frac{\ud x}{x}<\infty,
\]
then for all $\epsilon > 0$
\[
\p\Big(U(x)<(1-\epsilon)h(x),\textrm{ i.o., as }x\to
\infty\Big)=0.
\]
\end{itemize}
\end{corollary}
\section{The regular case}
The first type of tail behaviour of $I\big(\hat{\xi}\big)$ and
$\nu I\big(\hat{\xi}\big)$ that we consider is the case where
$\bar{F}$ and $\bar{F}_{\nu}$ satisfy
\begin{equation}\label{reg}
ct^{\alpha}L(t)\leq\bar{F}(t)\leq\bar{F}_{\nu}(t)\leq
Ct^{\alpha}L(t)\qquad\textrm{ as }\quad t\to 0,
\end{equation}
where $\alpha>0$, $c$ and $C$ are two positive constants such that
$c\leq C$ and $L$ is a slowly varying function at $0$. An important
example included in this case is when
$\bar{F}$ and $\bar{F}_{\nu}$ are regularly varying functions at $0$. \\
The ``regularity" of the behaviour  of $\bar{F}$ and
$\bar{F}_{\nu}$ gives us the following integral tests.
\begin{theorem}
Under condition (\ref{reg}), the lower envelope of $U$ at $0$ and
at $+\infty$ is as follows:
\begin{itemize}
\item[i)] Let $h\in\mathcal{H}^{-1}_{0}$, such that either
$\lim_{x\to 0}h(x)/x=0$ or $\liminf_{x\to 0}h(x)/x>0$, then
\[
\p\Big(U(x)<h(x), \textrm{ i.o., as } x\to 0\Big)= 0\textrm{ or }
1,
\]
according as
\[
\int_{0^+}\bar{F}\left(\frac{h(x)}{x}\right)\frac{\ud
x}{x}\qquad\textrm{is finite or infinite}.
\]\\
\item[ii)]Let $h\in\mathcal{H}^{-1}_{\infty}$, such that either
$\lim_{x\to +\infty}h(x)/x=0$ or $\liminf_{x\to +\infty}h(x)/x>0$,
then
\[
\p\Big(U(x)<h(x), \textrm{ i.o., as } x\to \infty\Big)= 0\textrm{
or } 1,
\]
according as
\[
\int^{+\infty}\bar{F}\left(\frac{h(x)}{x}\right)\frac{\ud
x}{x}\qquad\textrm{is finite or infinite}.
\]
\end{itemize}
\end{theorem}
{\it Proof:} First let us check that under condition (\ref{reg})
we  have
\begin{equation}\label{equivre0}
\int_{0}^{\lambda}\bar{F}_{\nu}\left(\frac{h(x)}{x}\right)\frac{\ud
x}{x}<\infty\quad\textrm{ if and only if
}\quad\int_{0}^{\lambda}\bar{F}\left(\frac{
h(x)}{x}\right)\frac{\ud x}{x}<\infty.
\end{equation}
Since $\nu I\leq I$ a.s., it is clear that we only need to prove
that
\[
\int_{0}^{\lambda}\bar{F}\left(\frac{h(x)}{x}\right)\frac{\ud
x}{x}<\infty\quad\textrm{ implies that
}\quad\int_{0}^{\lambda}\bar{F}_{\nu}\left(\frac{
h(x)}{x}\right)\frac{\ud x}{x}<\infty.
\]
From the hypothesis, either $\lim_{x\to 0}h(x)/x=0$ or
$\liminf_{x\to 0}h(x)/x>0$. In the first case, from condition
(\ref{reg}) there exists $\lambda>0$ such that, for every
$x<\lambda$
\[
c\left(\frac{h(x)}{x}\right)^{\alpha}L\left(\frac{h(x)}{x}\right)\leq
F\left(\frac{h(x)}{x}\right)\leq
C\left(\frac{h(x)}{x}\right)^{\alpha}L\left(\frac{h(x)}{x}\right).
\]
Since, we suppose that
$\int_{0}^{\lambda}F\left(\frac{h(x)}{x}\right)\frac{\ud x}{x}$ is
finite, then
\[
\int_{0}^{\lambda}\left(\frac{h(x)}{x}\right)^{\alpha}L\left(\frac{h(x)}{x}\right)\frac{\ud
x}{x}<\infty, \] and again from condition (\ref{reg}), we get that
$\int_{0}^{\lambda}F_{\nu}\left(\frac{h(x)}{x}\right)\frac{\ud
x}{x}$ is also finite. In the second case, since for any
$0<\delta<\infty$, $\p\big(I(\hat{\xi})<\delta\big)>0$, and
$\liminf_{x\to 0}h(x)/x>0$, we have for any $y$
 \begin{equation}\label{inprobreg}
0<\p\left(I\big(\hat{\xi}\big)<\liminf_{x\to
0}\frac{h(x)}{x}\right)<\p\left(
I\big(\hat{\xi}\big)<\frac{h(y)}{y}\right).
\end{equation}
Hence, since for every $t\geq 0$, $F(t)\leq F_{\nu}(t)$, we deduce
that
\[
\int_{0}^{\lambda}F\left(\frac{h(x)}{x}\right)\frac{\ud x}{x}=
\int_{0}^{\lambda}F_{\nu}\left(\frac{h(x)}{x}\right)\frac{\ud
x}{x}=\infty.
\]
Now, let us  check that for any constant $\beta>0$,
\begin{equation}\label{equivre}
\int_{0}^{\lambda}\bar{F}\left(\frac{h(x)}{x}\right)\frac{\ud
x}{x}<\infty\quad\textrm{ if and only if
}\quad\int_{0}^{\lambda}\bar{F}\left(\frac{\beta
h(x)}{x}\right)\frac{\ud x}{x}<\infty,
\end{equation}
Again, from the hypothesis either $\lim_{x\to 0}h(x)/x=0$ or
$\liminf_{x\to 0}h(x)/x>0$. In the first case, we deduce
(\ref{equivre}) from (\ref{reg}). In the
second case, from (\ref{inprobreg}) both of the integrals in (\ref{equivre}) are infinite.\\
Next, it follows from Theorem 3 part $(i)$ and (\ref{equivre0})
that if $\int_{0^+}\bar{F}\left(\frac{h(x)}{x}\right)\frac{\ud
x}{x}$ is finite, then for all $\epsilon>0$,
$\p\big(U(x)<(1-\epsilon)h(t), \textrm{ i.o., as }t\to 0\big)=0$.
If $\int_{0^+}\bar{F}\left(\frac{h(x)}{x}\right)\frac{\ud x}{x}$
diverges, then from Theorem 3 part $(ii)$ that for all
$\epsilon>0$, $\p\big(U(x)<(1+\epsilon)h(t), \textrm{ i.o., as
}t\to 0\big)=1$. Then $(\ref{equivre})$  allows
us to drop $\epsilon$ in this implications.\\
The tests at $+\infty$ are proven through the same way.\QED
Similarly, we get the following result for the upper envelope of
the future infimum.
\begin{theorem}
Under condition (\ref{reg}), the upper envelope of the future
infimum at $0$ and at $+\infty$ is as follows:
\begin{itemize}
\item[i)] Let $h\in\mathcal{H}_{0}$, such that either $\lim_{t\to
0}t/h(t)=0$ or $\liminf_{t\to 0}t/h(t)>0$, then
\[
\p\Big(J^{(0)}_{t}>h(t), \textrm{ i.o., as } t\to 0\Big)=
0\textrm{ or } 1,
\]
according as
\[
\int_{0^+}\bar{F}\left(\frac{t}{h(t)}\right)\frac{\ud
t}{t}<\infty\qquad\textrm{is finite or infinite}.
\]\\
\item[ii)]Let $h\in\mathcal{H}_{\infty}$, such that either
$\lim_{t\to +\infty}t/h(t)=0$ or $\liminf_{t\to +\infty}t/h(t)>0$,
then for all $x\geq 0$
\[
\p\Big(J_{t}^{(x)}>h(t), \textrm{ i.o., as } t\to \infty\Big)=
0\textrm{ or } 1,
\]
according as
\[
\int^{+\infty}\bar{F}\left(\frac{t}{h(t)}\right)\frac{\ud
t}{t}<\infty\qquad\textrm{is finite or infinite}.
\]
\end{itemize}
\end{theorem}
{\it Proof:} We prove this result by following the same arguments as
the proof of the previous Theorem.\QED An example of such a
behaviour will be given in section 7 (Example 3).
\section{The log regular case}
The second type of behaviour that we shall consider is when $\log
\bar{F}$ and $\log\bar{F}_{\nu}$ are regularly varying at $0$, i.e
\begin{equation}\label{logreg}
-\log\bar{F}_{\nu}(1/t)\sim-\log\bar{F}(1/t)\sim \lambda
t^{\beta}L(t),\quad \textrm{ as } t\to +\infty,
\end{equation}
where $\lambda>0$, $\beta>0$ and $L$ is a slowly varying function
at $+\infty$. Define the function
\begin{equation}\label{eqlogreg}
\psi(t)\eqdef\frac{t}{\inf\big\{s:1/\bar{F}(1/s)>|\log
t|\big\}},\quad t>0, \quad t\neq 1.
\end{equation}
Then the lower envelope of $U$ may be described as follows
\begin{theorem}Under condition (\ref{logreg}), the process $U$
satisfies the following law of the iterated logarithm:
\[
\liminf_{x\to 0}\frac{U(x)}{\psi(x)}=1,\quad \textrm{and}\quad
\liminf_{x\to +\infty}\frac{U(x)}{\psi(x)}=1\quad \textrm{a.s.}
\]
\end{theorem}
{\it Proof:} This Theorem is a consequence of Theorems 3 and 4,
and it is proven in the same way as Theorem 4 in \cite{ChP}, we
only need to emphasize that we can replace $\log \bar{F}_{\nu}$ by
$\log\bar{F}$, since they are asymptotically equivalent. \QED
Similarly the upper envelope of the future infimum may be
described as follows. Define the function
\begin{equation*}
\phi(t)\eqdef t\inf\big\{s:1/\bar{F}(1/s)>|\log t|\big\},\quad
t>0, \quad t\neq 1.
\end{equation*}
\begin{theorem}Under condition (\ref{logreg}), the future infimum process
satisfies the following law of the iterated logarithm:
\begin{itemize}
\item[i)]\[ \limsup_{t\to 0}\frac{J^{(0)}_{t}}{\phi(t)}=1,\qquad
\textrm{almost surely.}
\]
\item[ii)]For all $x\geq 0$,
\[ \limsup_{t\to
+\infty}\frac{J^{(x)}_{t}}{\phi(t)}=1,\qquad \textrm{almost
surely.}
\]
\end{itemize}
\end{theorem}
{\it Proof:} As the previous Theorem, the proof of this result
follows from Theorems 1 and 2, and we use exactly the same arguments
of the proof of Theorem 4 in \cite{ChP}.\QED In the following
section, we obtain sharp estimates for $\log \bar{F}$ and
$\log\bar{F}_{\nu}$ which will give us an important application of
this case.
\subsection{The case when $\xi$ has finite exponential moments}
Throughout this section we will suppose that the L\'evy process
$\xi$ associated to the PSSMP $X^{(x)}$ by its Lamperti
representation, has  finite exponential moments of arbitrary
positive order. This condition is satisfied, for example, when the
jumps of $\xi$ are bounded from above by some fixed number, in
particular when $\xi$ is a L\'evy process with no positive jumps.
Then, we have
\[
\e\big(e^{\lambda
\xi_{t}}\big)=\exp\big\{t\psi(\lambda)\big\}<+\infty\qquad t,
\lambda\geq 0.
\]
From Theorem 25.3 in Sato \cite{Sa}, we know that this hypothesis
is equivalent to assume that the L\'evy measure $\Pi$ of $\xi$
satisfies
\[ \int_{[1,\infty)}e^{\lambda x}\Pi(\ud x)
<+\infty\qquad \textrm{for every  } \lambda>0.
\]
Under this hypothesis, Bertoin and Yor \cite{BeY2} give a formula
for the negative moments of the exponential functional
$I(\hat{\xi})$
\begin{equation}\label{moments}
\e\Big(I\big(\hat{\xi}\big)^{-k}\Big)=m\frac{\psi(1)\cdots\psi(k-1)}{(k-1)!}\quad\textrm{
for }\quad k\geq 1,
\end{equation}
where $m=\e(\xi_{1})$ and with the convention that the right-hand
side equals $m$ for $k=1$. Moreover they proved that if $\xi$ has
no positive jumps, then $1/I(\hat{\xi})$ admits some exponential
moments, this means that the distribution of
$I(\hat{\xi})$ is determined by its negative entire moments.\\
From the entrance law of $X^{(x)}$ at $0$ (see (\ref{entrlaw})),
and the above equality (\ref{moments}),  we get the following
formula for the positive moments of $X^{(0)}_{1}$,
\begin{equation}\label{momentx}
\e_{0}\big(X_{1}^{k}\big)=\frac{\psi(1)\cdots\psi(k)}{k!}\quad
\textrm{for }\quad k\geq 1.
\end{equation}
Now, if we suppose that the Laplace exponent $\psi$ is regularly
varying at $+\infty$ with index $\beta \in (1,2)$, i.e.
$\psi(x)=x^{\beta}L(x)$, where $L$ is a slowly varying function at
$+\infty$; then from equation $(\ref{moments})$, we see
\[
\e\Big(I\big(\hat{\xi}\big)^{-k}\Big)=m\big((k-1)!\big)^{\beta
-1}L(1)\cdots L(k-1),
\]
and from $(\ref{momentx})$,
\[
\e_{0}\big(X_{1}^{k}\big)=\big(k!\big)^{\beta -1}L(1)\cdots L(k).
\]
In consequence, we can easily deduce that
\[
\e\bigg(\exp\Big\{\lambda
I^{-1}\big(\hat{\xi}\big)\Big\}\bigg)<+\infty\quad
\textrm{and}\quad \e_{0}\Big(\exp\{\lambda X_{1}\}\Big)<+\infty
\quad\textrm{ for all } \lambda>0.
\]
This allows us to apply the Kasahara's Tauberian Theorem (see
Theorem 4.12.7 in Bingham et al. \cite{Bing}) and get the
following estimate.
\begin{proposition}
Let $I(\hat{\xi})$ be the exponential functional associated to the
L\'evy process $\xi$. Suppose that $\psi$, the Laplace exponent of
$\xi$, varies regularly at $+\infty$ with index $\beta\in(1,2)$.
Then
\begin{equation}\label{tailp1}
-\log\p_{0}(X_{1}>x)\sim-\log
\p\left(I\big(\hat{\xi}\big)<1/x\right)\sim
(\beta-1)\overset{\leftharpoonup}{H}(x)\quad \textrm{ as }\quad
x\to +\infty,
\end{equation}
where
\begin{equation*}
\overset{\leftharpoonup}{H}(x)=\inf\Big\{s > 0\textrm{ , }
\psi(s)/s> x\Big\}.
\end{equation*}
\end{proposition}
Recall that if the process $\xi$ has no positive jumps then the fact
that the Laplace exponent $\psi$ is regularly varying at $\infty$
with index $\beta\in(1,2)$ is equivalent to that $\xi$ satisfies the
Spitzer's condition (see Proposition VII.6 in Bertoin~\cite{be}),
this is
\begin{equation*}
\lim_{t\to 0}\frac{1}{t}\int_{0}^{t}\p(\xi_{s}\geq 0)\ud
s=\frac{1}{\beta}.
\end{equation*}
\noindent\textit{Proof: } As we see above, the moment generating
functions of $I(\hat{\xi})^{-1}$ and  $X^{(0)}_{1}$ are well
defined for all $\lambda>0$. We will only prove the case of
$I(\hat{\xi})^{-1}$, the proof of the estimate of the tail
probability of $X^{(0)}_{1}$ is similar.\\
From the main result of
Geluk~\cite{Ge}, we know that if $\phi$ is a regularly varying
function at $+\infty$ with index $\sigma \in (0,1)$, then the
following are equivalent
\begin{itemize}
\item[(i)]
\hspace{1cm}$\bigg(\e\Big(I\big(\hat{\xi}\big)^{-n}\Big)/n!\bigg)^{1/n}\sim
e^{\sigma}/\phi(n)$ \hspace{.4cm} as  $n\to +\infty$,
\item[(ii)]\hspace{1cm} $\log\e\bigg(\exp\Big\{\lambda
I\big(\hat{\xi}\big)^{-1}\Big\}\bigg)\sim
\sigma\overset{\leftharpoonup}{\phi}(\lambda)$ \hspace{.4cm} as
$\lambda \to +\infty$,
\end{itemize}
where $\overset{\leftharpoonup}{\phi}(\lambda)=\inf\Big\{s >
0\textrm{ , } \phi(s)> \lambda\Big\}$.\\
If we have (ii), then a straightforward application of Kasahara's
Tauberian Theorem gives us
\begin{equation*}
-\log\p\Big(I\big(\hat{\xi}\big)^{-1}>x\Big)\sim
(1-\sigma)\overset{\leftharpoonup}{\varphi}(x)\qquad \textrm{as }
x\to +\infty,
\end{equation*}
where $\overset{\leftharpoonup}{\varphi}$ is the asymptotic inverse
of $s/\phi(s)$. Therefore, it is enough to show (i) with
$\phi(s)=s^{2}/\psi(s)$ to obtain the desired result.\\
Let us recall that if $\psi$ is regularly varying at $\infty$ with
index $\beta$, it can be expressed as $\psi(x)=x^{\beta}L(x)$,
where $L(x)$ is a slowly varying function. By the formula
(\ref{moments}) of negative moments of $I(\hat{\xi})$  and the
fact that $\psi$ is regularly varying, we have
\begin{equation*}
\bigg(\e\Big(I\big(\hat{\xi}\big)^{-n}\Big)/n!\bigg)^{1/n}=
m^{1/n}(n!)^{\frac{\beta-2}{n}}n^{\frac{1-\beta}{n}}\big(L(1)\cdots
L(n-1)\big)^{\frac{1}{n}},
\end{equation*}
due to the fact that $(n!)^{1/n}\sim ne^{-1}$ for $n$ sufficiently
large, then
\begin{equation*}
\bigg(\e\Big(I\big(\hat{\xi}\big)^{-n}\Big)/n!\bigg)^{1/n}\sim
\big(ne^{-1}\big)^{\beta-2}\exp\bigg\{\frac{1}{n}\sum_{k=
1}^{n}\log L(k)-\frac{1}{n}\log L(n)\bigg\}.
\end{equation*}
On the other hand, from the proof of Proposition 2 of Rivero
\cite{ri} we have that
\begin{equation*}
\frac{1}{n}\sum_{k= 1}^{n}\log L(k)\sim \log L(n) \qquad \textrm{
as} \quad n\to +\infty.
\end{equation*}
This implies that
\begin{equation*}
\bigg(\e\Big(I\big(\hat{\xi}\big)^{-n}\Big)/n!\bigg)^{1/n}\sim
e^{2-\beta}\frac{\psi(n)}{n^2}.
\end{equation*}
This last relation proves the Proposition.\QED Since the tail
probability of $I(\hat{\xi})^{-1}$ and $X^{(0)}_{1}$ have the same
asymptotic behaviour, it is logical to think that  the tail
probability of $(\nu I(\hat{\xi}))^{-1}$ could have the same
behaviour. The next Corollary confirms this last argument.
\begin{corollary}
Let $I(\hat{\xi})$ be the exponential functional associated to the
L\'evy process $\xi$. Suppose that $\psi$, the Laplace exponent of
$\xi$, varies regularly at $+\infty$ with index $\beta\in(1,2)$.
Then
\begin{equation*}
-\log \p\Big(\nu I\big(\hat{\xi}\big)<1/x\Big)\sim
(\beta-1)\overset{\leftharpoonup}{H}(x)\quad \textrm{ as }\quad
x\to +\infty,
\end{equation*}
where
\begin{equation*}
\overset{\leftharpoonup}{H}(x)=\inf\Big\{s > 0\textrm{ , }
\psi(s)/s> x\Big\}.
\end{equation*}
\end{corollary}
\textit{Proof: } Since $\nu I(\hat{\xi})\leq I(\hat{\xi})$ a.s.,
then
\[
-\log\p\Big(\nu I\big(\hat{\xi}\big)<1/x\Big)\leq
-\log\p\Big(I\big(\hat{\xi}\big)<1/x\Big).
\]
On the other hand, from the scaling property and since
$X^{(0)}_{1}\geq J^{(0)}_{1}$ a.s., we see
\[ -\log\p\Big(\nu
I\big(\hat{\xi}\big) <1/x\Big)=-\log\p\big( U(1)<1/x\big)\geq -\log
\p_{0}(X_{1}>x).
\]
Hence, from the estimate (\ref{tailp1}) we have that
\[
-\log\p\Big(\nu I\big(\hat{\xi}\big) <1/x\Big)\sim
(\beta-1)\overset{\leftharpoonup}{H}(x)\quad \textrm{ as }\quad
x\to +\infty,
\]
and this finishes the proof.\QED These estimates will allow us to
obtain laws of iterated logarithm for the last
passage time process and for the future infimum process in terms of the following function.\\
Let us define the function
\[
h(t)=\frac{\log|\log t|}{\psi(\log|\log t|)}\qquad
\textrm{for}\qquad t>1,\quad t\neq e.
\]
By integration by parts, we can see that the function
$\psi(\lambda)/\lambda$ is increasing, hence it is straightforward
that the function $th(t)$ is also increasing in a neighbourhood of
$\infty$.
\begin{corollary}If $\psi$ is regularly varying at $+\infty$ with index $\beta\in
(1,2)$, then
\begin{equation}\label{lawith1}
\liminf_{x\to 0}\frac{U(x)}{xh(x)}=(\beta-1)^{\beta -1} \qquad
\textrm{almost surely}
\end{equation}
and,
\begin{equation}\label{lawith3}
\liminf_{x\to +\infty}\frac{U(x)}{xh(x)}=(\beta-1)^{\beta -1} \qquad
\textrm{almost surely.}
\end{equation}
\end{corollary}
\textit{Proof:} It is enough to see that for $t$ sufficiently
small and $t$ sufficiently large the functions $th(t)$ and $\psi
(t)$ are asymptotically equivalent, but this is clear from
(\ref{logreg}). Now, applying Theorem 7 we obtain the desired law
of the iterated logarithm.\QED Let us define
\[
f(t)=\frac{\psi(\log |\log t|)}{\log |\log t|}\qquad{\rm for}\quad
t>1, \quad t\neq e.
\]
\begin{corollary}If $\psi$ is regularly varying at $+\infty$ with index $\beta\in
(1,2)$, then
\begin{equation}\label{lawith12}
\limsup_{t\to 0}\frac{J_{t}}{tf(t)}=(\beta-1)^{-(\beta -1)} \qquad
\textrm{almost surely}
\end{equation}
and,
\begin{equation}\label{lawith32}
\limsup_{t\to +\infty}\frac{J(t)}{tf(t)}=(\beta-1)^{-(\beta -1)}
\qquad \textrm{almost surely}.
\end{equation}
\end{corollary}
{\it Proof:} This proof follows from similar arguments of the last
Corollary and using Theorem 8.\QED\noindent{\bf Example 1} Let
$X^{(0)}_{t}$ be a stable L\'evy process conditioned to stay
positive with no positive jumps and with index $1<\alpha\leq 2$,
(see Bertoin \cite{be} for a proper definition). From Theorem VII.18
in \cite{be}, we know that the process time-reversed at its last
passage time below $x,$ $(x-X^{(0)}_{(U(x)-t)^{-}}, 0\leq t\leq
U(x))$, has the same law as the killed process at its first passage
time above $x$, $(\xi_{t}, 0\leq t\leq T_{x})$, where $\xi$ is a
stable L\'evy process with no
positive jumps and with the same index as $X^{(0)}$.\\
From Theorem VII.1 in \cite{be}, we know that $(T_{x},x\geq 0)$ is a
subordinator with Laplace exponent
$\Phi(\lambda)=\lambda^{1/\alpha}$. Hence by the previous argument,
we will have that $X^{(0)}$ drifts towards $+\infty$ and that the
process $(U(x), x\geq 0)$ is a stable subordinator with index
$1/\alpha$. Hence an application of the Tauberian theorem of de
Brujin (see for instance Theorem 5.12.9 in Bingham et al.
\cite{Bing}) gives us the following estimate
\[
-\log \bar{F}(x)\sim
\frac{\alpha-1}{\alpha}\left(\frac{1}{\alpha}\right)^{1/(\alpha-1)}x^{-1/(\alpha
-1)}\quad\textrm{ as }\quad x\to 0.
\]
Note that due to the absence of positive jumps $\nu=1$ a.s.\\
Then applying Theorems 7 and 8, we get the following law of the
iterated logarithm.
\begin{corollary}
Let $X^{(0)}$ be a stable L\'evy process conditioned to stay
positive  with no positive jumps and $\alpha>1$. Then, its related
last passage time process satisfies
\[
\liminf_{x\to 0}\frac{U(x)\big(\log|\log
x|\big)^{\alpha-1}}{x^{\alpha}}=\frac{1}{\alpha}\left(1-\frac{1}{\alpha}\right)^{\alpha-1},\qquad\textrm{almost
surely}.
\]
The same law of the iterated logarithm is satisfied for large
times.\\ The future infimum process of $X^{(0)}$ also satisfies
\[
\limsup_{t\to 0}\frac{J^{(0)}_{t}}{t^{1/\alpha}\big(\log|\log
x|\big)^{1-1/\alpha}}=\alpha(\alpha-1)^{-\frac{\alpha-1}{\alpha}},\quad\textrm{
almost surely},
\]
and for all $x\geq 0$,
\[
\limsup_{t\to +\infty}\frac{J^{(x)}_{t}}{t^{1/\alpha}\big(\log|\log
x|\big)^{1-1/\alpha}}=\alpha(\alpha-1)^{-\frac{\alpha-1}{\alpha}},\quad\textrm{
almost surely}.
\]
\end{corollary}
{\bf Example 2} Let us
suppose that $\xi=(Y_{t}+ct, t\geq 0)$, where $Y$ is a stable
L\'evy process of index $\beta\in (1,2)$ with no positive jumps
and $c$ a positive constant. Its Laplace exponent has the form
\[
\e\big(e^{\lambda
\xi_{t}}\big)=\exp\{t(\lambda^{\beta}+c\lambda)\}, \qquad
\textrm{for }t\geq 0,\textrm{ and } \lambda> 0,
\]
where $c>0$. Note that under  the hypothesis that $Y$ has no
positive jumps, $\nu=1$ a.s. Let us define by $X^{(x)}$, the PSSMP
associated to $\xi$ starting from $x$ and with scaling index
$\alpha>0$, then when $x=0$ its last passage process $U$ satisfies
\[
\liminf_{x\to 0}\frac{U(x)}{x^{\alpha}\big(\log|\log
x|\big)^{(1-\beta)\alpha }}=\alpha^{-\beta\alpha}(\beta
-1)^{\alpha(\beta-1)},\quad \textrm{almost surely.}
\]
We have the same law of the iterated logarithm at $+\infty$.\\
The future infimum process $J^{(x)}$ satisfies that
\[
\limsup_{t\to 0
}\frac{J^{(0)}_{t}}{t^{\frac{1}{\alpha}}\big(\log|\log
t|\big)^{\frac{(\beta-1)}{\alpha}
}}=\alpha^{\frac{\beta}{\alpha}}(\beta
-1)^{-\frac{\beta-1}{\alpha}},\quad \textrm{almost surely,}
\]
and for all $x\geq 0$
\[
\limsup_{t\to
+\infty}\frac{J^{(x)}_{t}}{t^{\frac{1}{\alpha}}\big(\log|\log
t|\big)^{\frac{\beta-1}{\alpha}}}=\alpha^{\frac{\beta}{\alpha}}(\beta
-1)^{-\frac{\beta-1}{\alpha}},\quad \textrm{almost surely.}
\]
Note that when $\alpha=\beta$, the process $J$ has the same
asymptotic behaviour as $\xi$, this is
\[
\limsup_{t\to 0 (\textrm{or }
\mbox{}+\infty)}\frac{\xi_{t}}{t^{\frac{1}{\beta}}\big(\log|\log
t|\big)^{\frac{\beta-1}{\beta}}}=\beta(\beta-1)^{-\frac{\beta-1}{\beta}},\quad
\textrm{almost surely},
\]
see Zolotarev \cite{zo} for details, and also the same asymptotic
behaviour of $\bar{Y}$, the stable L\'evy process conditioned to
stay positive with no positive jumps (see Corollary 3).
\section{Transient Bessel
process.} In this section we will suppose that $\xi=(2(B_{t}+at),
t\ge 0)$,
where $B$ is a standard Brownian motion and $a>0$.\\
We define the process $Z=(Z_{t}, t\geq 0)$, the square of the
$\delta$-dimensional Bessel processes starting at $x\geq 0$, as the
unique strong solution of the stochastic differential equation
\begin{equation}\label{SDE}
Z_{t}=x+2\int_{0}^{t}\sqrt{|Z_{s}|}\ud \beta_{s}+\delta
t,\qquad\textrm{for} \quad \delta\geq 0,
\end{equation}
where $\beta$ is a standard Brownian motion.\\
By the Lamperti representation, we know that we can define $X^{(x)}$
a PSSMP starting at $x>0$, such that
\[
X^{(x)}_{xI_{t}(\xi)}=x\exp\big\{\xi_{t}\big\}\quad\textrm{for}\quad
t\geq 0.
\]
Then, applying the It\^o's formula and Dubins-Schwartz's Theorem
(see for instance, Revuz and Yor~\cite{RY}), we get
\[
X^{(x)}_{xI_{t}(\xi)}=x+2\int_{0}^{xI_{t}(\xi)}\sqrt{X^{(x)}_{s}}\ud
B_{s}+2(a+1)xI_{t}(\xi).
\]
Hence it follows that $X^{(x)}$ satisfies (\ref{SDE}) with
$\delta=2(a+1)$ and therefore $X^{(x)}$ is the square of the
$\delta$-dimensional Bessel processes starting at $x>0$. From the
main result of Bertoin and Yor ~\cite{BeY}, we may define $X^{(x)}$
at $x=0$, and we can computed its entrance law by (\ref{entrlaw}).
Since, we
suppose that $a>0$, we deduce that $X^{(x)}$ is a transient process and that $\delta>2$.\\
From the formula of negative moments (\ref{moments}) of the
exponential functional $I(\hat{\xi})$, we can deduce (see Example
3 in Bertoin and Yor ~\cite{BeY2}) the following identity in
distribution
\begin{equation}\label{igbessel}
\int_{0}^{\infty}\exp\big\{-2(B_{s}+as)\big\}\ud s \ed
\frac{1}{2\gamma_{a}},
\end{equation}
where $\gamma_{a}$ is a gamma random variable with index $a>0$. In
fact, we can also deduce that $X^{(0)}_{1}$
is distributed as $2\gamma_{a+1}$. \\
We recall that the distribution of $\gamma_{a}$ for $a>0$, is given
by
\begin{equation}\label{igbessel1}
\p (\gamma_{a}\leq x)=\frac{1}{\Gamma(a)}\int_{0}^{x}e^{-y}y^{a
-1}\ud y,\quad \textrm{where
}\quad\Gamma(a)=\int_{0}^{\infty}e^{-y}y^{a -1}\ud y.
\end{equation}
It is important to note that due the continuity of the paths of
$X^{(0)}$, we have that $\nu=1$ almost surely.\\
The following Lemma will be helpful for the application of our
general results to the case of transient Bessel processes.
\begin{lemma}
Let $a>0$, then there exist $c$ and $C$, two positive constants such
that
\[
ce^{-x}x^{a-1}\leq\int_{x}^{\infty}e^{-y}y^{a -1}\ud y\leq
Ce^{-x}x^{a-1}, \quad \textrm{for} \quad x\geq \frac{C(a-1)}{C-1}.
\]
\end{lemma}
From this Lemma, we deduce the following integral tests.
\begin{theorem}
Let $h\in \mathcal{H}^{-1}_{0}$, then:
\begin{itemize}
\item[i)] If
\[
\int_{0^{+}}\Big(x/2h(x)
\Big)^{\frac{\delta-4}{2}}\exp\Big\{-x/2h(x)\Big\}\frac{\ud
x}{x}<\infty,
\]
then for all $\epsilon > 0$
\[
\p\Big(U(x)<(1-\epsilon)h(x),\textrm{ i.o., as }x\to 0\Big)=0.
\]
\item[ii)] If
\[
\int_{0^{+}}\Big(x/2h(x)
\Big)^{\frac{\delta-4}{2}}\exp\Big\{-x/2h(x)\Big\}\frac{\ud
x}{x}=\infty,
\]
then for all $\epsilon > 0$
\[
\p\Big(U(x)<(1+\epsilon)h(x),\textrm{ i.o., as }x\to 0\Big)=1.
\]
\end{itemize}
\end{theorem}
{\it Proof}: The proof of this Theorem follows from the fact that
\begin{equation}\label{iggama}
\p(I<x)=\p\big(\gamma_{(\delta-2)/2}>1/2x\big) \quad \textrm{for
}\quad x>0,
\end{equation}
and an application of Theorem 3 and Lemma 2.\QED
\begin{theorem}
Let $h\in \mathcal{H}^{-1}_{\infty}$, then:
\begin{itemize}
\item[i)] If
\[
\int^{+\infty}\Big(x/2h(x)
\Big)^{\frac{\delta-4}{2}}\exp\Big\{-x/2h(x)\Big\}\frac{\ud
x}{x}<\infty,
\]
then for all $\epsilon > 0$
\[
\p\Big(U(x)<(1-\epsilon)h(x),\textrm{ i.o., as }x\to
+\infty\Big)=0.
\]
\item[ii)] If
\[
\int^{+\infty}\Big(x/2h(x)
\Big)^{\frac{\delta-4}{2}}\exp\Big\{-x/2h(x)\Big\}\frac{\ud
x}{x}=\infty,
\]
then for all $\epsilon > 0$
\[
\p\Big(U(x)<(1+\epsilon)h(x),\textrm{ i.o., as }x\to
+\infty\Big)=1.
\]
\end{itemize}
\end{theorem}
{\it Proof:} The proof of these integral tests is very similar to
the proof of the previous result, it is enough to apply Lemma 2 and
Theorem 4 to equality (\ref{iggama}).\QED From these integral tests,
we get the following law of iterated logarithm.
\[
\liminf_{x\to 0}U(x) \frac{2\log|\log x|}{x}=1\quad\textrm{and}\quad
\liminf_{x\to +\infty}U(x) \frac{2\log\log x}{x}=1\quad\textrm{
almost surely.}
\]
Note that we are also in the ``logregular" case and we can apply
Theorem 7  to get the same law of the iterated logarithm.\\
For the upper envelope of the future infimum process, we have the
following integral tests.
\begin{theorem}
Let $h\in \mathcal{H}_0$, then:
\begin{itemize}
\item[i)] If
\[
\int_{0^{+}}\Big(h(t)/2t
\Big)^{\frac{\delta-4}{2}}\exp\Big\{-h(t)/2t\Big\}\frac{\ud
t}{t}<\infty,
\]
then for all $\epsilon > 0$
\[
\p\Big(J^{(0)}_t>(1+\epsilon)h(t),\textrm{ i.o., as }t\to
0\Big)=0.
\]
\item[ii)] If
\[
\int_{0^{+}}\Big(h(t)/2t
\Big)^{\frac{\delta-4}{2}}\exp\Big\{-h(t)/2t\Big\}\frac{\ud
t}{t}=\infty,
\]
then for all $\epsilon > 0$
\[
\p\Big(J^{(0)}_t>(1-\epsilon)h(t),\textrm{ i.o., as }t\to
0\Big)=1.
\]
\end{itemize}
\end{theorem}
{\it Proof}: We get this result applying Theorem 1 and Lemma 2 to
the equality (\ref{iggama}).\QED
\begin{theorem}
Let $h\in \mathcal{H}_{\infty}$, then for all $x\geq 0$:
\begin{itemize}
\item[i)] If
\[
\int^{+\infty}\Big(h(t)/2t
\Big)^{\frac{\delta-4}{2}}\exp\Big\{-h(t)/2t\Big\}\frac{\ud
t}{t}<\infty,
\]
then for all $\epsilon > 0$
\[
\p\Big(J^{(x)}_t>(1+\epsilon)h(t),\textrm{ i.o., as }t\to
+\infty\Big)=0.
\]
\item[ii)] If
\[
\int^{+\infty}\Big(h(t)/2t
\Big)^{\frac{\delta-4}{2}}\exp\Big\{-h(t)/2t\Big\}\frac{\ud
t}{t}=\infty,
\]
then for all $\epsilon > 0$
\[
\p\Big(J^{(x)}_t>(1-\epsilon)h(t),\textrm{ i.o., as }t\to
+\infty\Big)=1.
\]
\end{itemize}
\end{theorem}
{\it Proof: } The proof of these integral test is similar to the
previous Theorem. We only replace Theorem 1 by Theorem 2. \QED
From these integral tests, we get the following law of iterated
logarithm,
\[
\limsup_{t\to 0}\frac{J^{(0)}_t}{2t\log|\log
t|}=1\quad\textrm{and}\quad \limsup_{t\to
+\infty}\frac{J^{(x)}_t}{2t\log\log t}=1\quad\textrm{ almost
surely,}
\]
for $x\geq 0$.  Here we can also obtain the same law of the iterated
logarithm applying Theorem 8.
\section{The upper envelope of increasing self-similar Markov processes.}
From the Lamperti representation (\ref{lamp}), we know that
increasing PSSMP are related to subordinators and by the main result
of Bertoin and Caballero we can define this processes at $0$. In
this
section we suppose that $\xi$ is a subordinator.\\
Since the PSSMP $X^{(x)}$ is increasing, we know that its supremum ,
its past infimum and its future infimum are the same. Moreover, its
first passage time over the level $y>0$ is the same as the last
passage time below $y$. Therefore, with our previous main results we
can describe the upper envelope of $X^{(x)}$ at $0$ (when $x=0$) and
at $+\infty$ (for all $x\geq 0$) and also the lower envelope of its
first passage time defined by,
\[
S_y=\inf\Big\{t: X^{(0)}_t\geq y\Big\},\quad\textrm{ for } y>0.
\]
From Proposition 2.1 in Carmona, Petit and Yor \cite{CPY} provided
that $m<\infty$, we know that the law of $I(\hat{\xi})$ admits a
density $\rho$ which is infinitely differentiable on
$(0,+\infty)$. Moreover, from (\ref{entrlaw}) the entrance law has
also a density and is described as follows
\begin{equation}\label{density}
\rho_{1}(x)=\frac{1}{m x}\rho\left(\frac{1}{x}\right)\qquad
\textrm{for } x\in(0, \infty).
\end{equation}
Generally speaking, we cannot estimate the tail probability
$\p(I(\hat{\xi})<t)$ for $t$ sufficiently small, so we will now
give some examples in which we can obtain an estimation.
\subsection{Examples}
{\bf Example 1} Let $\xi=N$ be a standard Poisson process. From
Proposition 3 in Bertoin and Yor \cite{BeY2}, we know that
\[
-\log\p\Big(I\big(\hat{\xi}\big)<t\big)\sim\frac12(\log1/t)^2,
\qquad \textrm{as }t\to 0,
\]
and also that
\[
-\log\rho(x)\sim\frac12(\log1/x)^2, \qquad \textrm{as }x\to 0.
\]
From (\ref{density}) we get that
\[
-\log\rho_{1}(x)\sim\frac12(\log x)^2, \qquad \textrm{as }x\to
+\infty.
\]
Now, applying Theorem 4.12.10 in Bingham et al. \cite{Bing} and
doing some elementary calculations,  we obtain that
\[
-\log \int_{x}^{+\infty}\rho_{1}(y)\ud y \sim\frac12(\log
x)^2\qquad \textrm{ as } x\to +\infty.
\]
These estimations allow us the following laws of the iterated
logarithm. Let us define
\[
f(t)=t\exp\Big\{-\sqrt{2\log|\log t|}\Big\}.
\]
\begin{corollary} Let $N$ be a Poisson process, then the PSSMP $X^{(x)}$
associated to $N$ by the Lamperti representation satisfies the
following law of the iterated logarithm,
\[
\limsup_{t\to 0}\frac{X^{(0)}_{t}f(t)}{t^2}=1, \quad\textrm{
almost surely.}
\]
For all $x\geq 0$,
\[
\limsup_{t\to +\infty}\frac{X^{(x)}_{t}f(t)}{t^2}=1, \quad\textrm{
almost surely.}
\]
The first passage time process $S$ associated to $X^{(0)}$
satisfies the following law of the iterated logarithm,
\[
\liminf_{x\to 0}\frac{S_x}{f(x)}=1, \quad\textrm{and}\quad
\liminf_{x\to +\infty}\frac{S_x}{f(x)}=1 \quad \textrm{ a. s.}
\]

\end{corollary}
{\it Proof:} The proof of this Corollary is consequence of the
previous estimations and application of the divergent parts of
Theorems 1, 2, 3 and 4, and Corollaries 1 and 2. \QED \noindent{\bf
Example 2 }Let $\xi$ be a subordinator with zero drift and L\'evy
measure $\Pi(\ud x)=abe^{-bx}\ud x$, with $a,b>0$, i.e. a compound
Poisson process with jumps having an exponential distribution.
Carmona, Petit and Yor showed that the density $\rho$ of
$I\big(\hat{\xi}\big)$ is given by
\[
\rho(x)=\frac{a^{1+b}}{\Gamma(1+b)}x^be^{-ax},\qquad\textrm{ for
}x>0.
\]
The PSSMP associated to $\xi$ by the Lamperti representation is the
well-know generalized Watanabe process. From (\ref{density}), we get
that
\[
\p_{0}\big(X_1>y\big)=\frac{ba^{1+b}}{\Gamma(1+b)}\int_0^{1/y}z^{b-1}e^{-az}
\ud z.
\]
On the other hand, It is clear that
\[
\p\Big(I\big(\xi\big)<y\Big)= \frac{a^{1+b}}{\Gamma(1+b)}\int_0^y
x^b e^{-ax}\ud x
\]
Elementary calculations give us the following inequality,
\[
e^{-ax}\frac{x^{b+1}}{b+1}\leq \int_0^x z^b e^{-az}\ud z\leq
x^b\frac{(1-e^{-ax})}{a}.
\]
Hence for $x$ sufficiently small, there exists $c_b$ a positive
constant such that
\[
\p\Big(I\big(\xi\big)<x\Big)\sim c_b
\frac{a^{1+b}}{\Gamma(1+b)}x^{b+1}e^{-ax},
\]
and for $y$ sufficiently large there exist $C_b$ such that
\[
\p_0\Big(X_1>y\Big)\sim C_b
\frac{a^{1+b}}{\Gamma(1+b)}(1/y)^{b}e^{-a/y}.
\]
Then applying the divergent part of Theorems 1, 2, 3 and 4 and
Corollaries 1 and 2, we get the following laws of the iterated
logarithm for the generalized Watanabe process and its first
passage time. Let us define
\[
g(t)=a^{-1}t\log|\log t|.
\]
\begin{corollary}
Let $\xi$  be a compound Poisson process with jumps having and
exponential distribution as above, then the PSSMP $X^{(x)}$
associated to $\xi$ by the Lamperti representation satisfies the
following law of the iterated logarithm,
\[
\limsup_{t\to 0}\frac{X^{(0)}_{t}g(t)}{t^2}=1, \quad\textrm{
almost surely.}
\]
For all $x\geq 0$,
\[
\limsup_{t\to +\infty}\frac{X^{(x)}_{t}g(t)}{t^2}=1, \quad\textrm{
almost surely.}
\]
The first passage time process $S$ associated to $X^{(0)}$
satisfies the following law of the iterated logarithm,
\[
\liminf_{x\to 0}\frac{S_x}{g(x)}=1, \quad\textrm{and}\quad
\liminf_{x\to 0}\frac{S_x}{g(x)}=1 \quad\textrm{ a. s.}
\]
\end{corollary}
\noindent{\bf Example 3} Let $\xi$ a subordinator with zero drift
and L\'evy measure
\[
\Pi(\ud x)=\frac{\beta e^x}{\Gamma(1-\beta)(e^x-1)^{1+\beta}}\ud
x,
\]
with $\beta\in(0,1)$. The PSSMP $X^{(x)}$ associated to $\xi$ is
the stable subordinator of index $\beta$ (see for instance Rivero \cite{ri}). \\
From Zolotarev \cite{zo1}, we know that there exists $k$ a positive
constant such that
\[
\p_0(X_1>x)\sim kx^{-\beta} \qquad x\to +\infty.
\]
It is well-known that the law of $X^{(0)}_1$ has a density $\rho_1$
with respect to the Lebesgue measure and that this density is
unimodal, i.e., there exist $b>0$ such that $\rho_1(x)$ is
increasing in $(0, b)$ and decreasing in $(b,+\infty)$ (see for
instance Sato \cite{Sa}). Hence $\rho_1$ is monotone in a
neighborhood of $+\infty$, then by the monotone density Theorem (see
Theorem 1.7.2 in Bingham et al.\cite{Bing} page 38) we get
\[
\rho_1(x)\sim k\beta x^{-\beta-1}\qquad x\to +\infty.
\]
From (\ref{density}), we can easily deduce that
\[
\rho(x)\sim m k\beta x^{\beta}\qquad x\to 0,
\]
and it is also easy to see that
\[
\p\Big(I\big(\hat{\xi}\big)<x\Big)\sim m k\beta x^{\beta+1}\qquad
x\to 0.
\]
Hence, we can apply Theorems 5 and 6, and deduce the following
results.
\begin{corollary}
Let $\xi$ be a subordinator without drift and such that its L\'evy
mesure $\Pi$ satisfies
\[
\Pi(\ud x)=\frac{\beta e^x}{\Gamma(1-\beta)(e^x-1)^{1+\beta}}\ud x.
\]
The lower envelope of $S$, the first passage time of the PSSMP
$X^{(0)}$, at $0$ and at $+\infty$ is as follows:
\begin{itemize}
\item[i)] Let $h\in\mathcal{H}^{-1}_{0}$, such that either
$\lim_{x\to 0}h(x)/x=0$ or $\liminf_{x\to 0}h(x)/x>0$, then
\[
\p\Big(S_x<h(x), \textrm{ i.o., as } x\to 0\Big)= 0\textrm{ or }
1,
\]
according as
\[
\int_{0^+}\left(\frac{h(x)}{x}\right)^{\beta+1}\frac{\ud x}{x}\quad
\textrm{ is finite or infinite}.
\]\\
\item[ii)]Let $h\in\mathcal{H}^{-1}_{\infty}$, such that either
$\lim_{x\to +\infty}h(x)/x=0$ or $\liminf_{x\to +\infty}h(x)/x>0$,
then
\[
\p\Big(S_x<h(x), \textrm{ i.o., as } x\to \infty\Big)= 0\textrm{
or } 1,
\]
according as
\[
\int^{+\infty}\left(\frac{h(x)}{x}\right)^{\beta+1}\frac{\ud
x}{x}\quad \textrm{ is finite or infinite}.
\]
\end{itemize}
\end{corollary}
\begin{corollary}
Let $\xi$ be a subordinator as in Corollary 8,  the upper envelope
of $X^{(x)}$ at $0$ ($x=0$) and at $+\infty$ ($x\geq 0$) is as
follows:
\begin{itemize}
\item[i)] Let $h\in\mathcal{H}_{0}$, such that either $\lim_{t\to
0}t/h(t)=0$ or $\liminf_{t\to 0}t/h(t)>0$, then
\[
\p\Big(X^{(0)}_{t}>h(t), \textrm{ i.o., as } t\to 0\Big)=
0\textrm{ or } 1,
\]
according as
\[
\int_{0^+}\left(\frac{x}{h(x)}\right)^{\beta+1}\frac{\ud x}{x} \quad
\textrm{ is finite or infinite}.
\]
\\
\item[ii)]Let $h\in\mathcal{H}_{\infty}$, such that either
$\lim_{t\to +\infty}t/h(t)=0$ or $\liminf_{t\to +\infty}t/h(t)>0$,
then for all $x\geq 0$
\[
\p\Big(X_{t}^{(x)}>h(t), \textrm{ i.o., as } t\to \infty\Big)=
0\textrm{ or } 1,
\]
according as
\[
\int^{+\infty}\left(\frac{x}{h(x)}\right)^{\beta+1}\frac{\ud
x}{x}\quad \textrm{ is finite or infinite}.
\]
\end{itemize}
\end{corollary}
{\bf Acknowledgements.} I would like to thank Lo\"\i c Chaumont
for guiding me through the development of this work, and for all
his helpful advice. I would also like to thank Victor Rivero for
all his useful suggestions.

\end{document}